\documentclass{article}
\usepackage{amsmath}
\usepackage{graphicx}
\usepackage{subfig}
\usepackage{color}
\usepackage{placeins}
\usepackage{authblk}
\usepackage{tikz}
\graphicspath{ {./pics/} }

\newcommand{\da}{\textcolor{black}}

\title{Multicontinuum Modeling of Time-Fractional Diffusion-Wave Equation in Heterogeneous Media}

\date{}

\author[a]{Huiran Bai}
\author[b]{Dmitry Ammosov}
\author[a]{Yin Yang \thanks{Corresponding author}}
\author[a]{Wei Xie}
\author[b]{Mohammed Al Kobaisi}

\affil[a]{School of Mathematics and Computational Science, Xiangtan University, National Center for Applied Mathematics in Hunan, Xiangtan 411105, Hunan, China}
\affil[b]{Chemical and Petroleum Engineering Department, Khalifa University of Science and Technology, Abu Dhabi, 127788, UAE}

\begin{document}
	\bibliographystyle{unsrt}
	\maketitle






\renewcommand{\thefootnote}{}%
\footnotetext{E-mail addresses: 202331510101@smail.xtu.edu.cn (Huiran Bai), dmitrii.ammosov@ku.ac.ae (Dmitry Ammosov), yangyinxtu@xtu.edu.cn (Yin Yang), xiew@smail.xtu.edu.cn (Wei Xie), mohammed.alkobaisi@ku.ac.ae (Mohammed Al Kobaisi).}
\addtocounter{footnote}{-1}%



\begin{center}
	\textbf{Abstract}
\end{center} 

This paper considers a time-fractional diffusion-wave equation with a high-contrast heterogeneous diffusion coefficient. A numerical solution to this problem can present great computational challenges due to its multiscale nature. Therefore, in this paper, we derive a multicontinuum time-fractional diffusion-wave model using the multicontinuum homogenization method. For this purpose, we formulate constraint cell problems considering various homogenized effects. These cell problems are implemented in oversampled regions to avoid boundary effects. By solving the cell problems, we obtain multicontinuum expansions of fine-scale solutions. Then, using these multicontinuum expansions and supposing the smoothness of the macroscopic variables, we rigorously derive the corresponding multicontinuum model. Finally, we present numerical results for two-dimensional model problems with different time-fractional derivatives to verify the accuracy of our proposed approach.


	\section{Introduction}
	
	
	Many chemical and physical processes involve memory \da{and} genetic effects. \da{It is known that} fractional differential equations are better suited to modeling such processes than integer differential equations \cite{zheng2020legendre,2sun2016some}. \da{Time-}fractional diffusion-wave equations are generalizations of classical parabolic $(\alpha=1)$ \da{and} hyperbolic equation\da{s} $(\alpha=2)$, which are derived \da{by} replacing the first- or second-order time derivative \da{with} a fractional derivative of order $\alpha$ with $0<\alpha<2$ \cite{mainardi1997fractional}. In the case $0<\alpha<1$, \da{we have} \da{the time-}fractional diffusion or subdiffusion equation, while in the case $1<\alpha<2$, \da{we obtain the} superdiffusion equation, \da{or the time-}fractional wave equation\da{,} or \da{the time-}fractional diffusion-wave equation. The time-fractional diffusion-wave equation describes important physical phenomena that occur in amorphous, colloid, glass and porous materials, fractal and percolation clusters, comb structures, dielectric and semiconductor, biological systems, polymers, random and disordered media, geophysical and geological processes and the universal electromagnetic, acoustic, and mechanical \da{processes} \da{(see} \cite{metzler2000random,gafiychuk2006pattern,hilfer2000applications,metzler2004restaurant,nigmatullin1984theoretical,nigmatullin1986realization}). 
	
	In the past two decades, \da{there have been many works dedicated to} the time-fractional diffusion-wave equation (TFDWE). Mainardi and Paradisi \cite{mainardi2001fractional} applied TFDWE to study the propagation of stress waves in viscoelastic media relevant to acoustics and seismology. Sun and Wu \cite{sun2006fully} proposed a fully discrete difference scheme for TFDWE by introducing two new variables to transform the original equation into a low\da{-}order system of equations. Du et al. \cite{du2010compact} proposed a compact difference scheme for the fractional diffusion-wave equation. Zhang et al. \cite{zhang2012compact} presented the alternating direction implicit (ADI) method in the time stepping and a difference scheme combining the compact difference approach for spatial discretization for solving the two-dimensional TFDWE. Li and Li \cite{li2021fast} presented a fast element-free Galerkin (EFG) method. Alikhanov \cite{alikhanov2024second} presented the investigation of a TFDWE involving a Caputo derivative and examined a more generalized model incorporating a memory kernel and the generalized Riemann–Liouville fractional integral.


	However, many real problems have high contrast properties and complex heterogeneities at multiple scales. \da{Examples include} composite materials, porous media \cite{brown2016generalized}, turbulent transport in high Reynolds number flows \cite{baldwin1991one}, {and others.} Adopting traditional numerical methods for \da{such multiscale problems} is \da{challenging} even \da{using} supercomputers. Simulations of these problems are \da{costly} and can cause difficulties \da{due to} a large number of degrees of freedom. \da{To reduce the computational cost, various upscaling techniques are often used to perform simulations on coarse computational grids.}

	
	Many different \da{upscaling} approach\da{e}s have been exploited to \da{solve} multiscale problems on a coarse computational grid. \da{Among these approaches, one can mention} homogenization methods \cite{wu2002analysis,hornung2012homogenization}. The main idea of the homogenization methods is to compute homogenized coefficients \da{(effective properties) by solving} cell problems and derive macroscopic equations on \da{a} coarse grid. \da{This significantly reduces computational costs and speeds up calculations. However, in the case of high-contrast heterogeneous media, the classical homogenization method may not be sufficient} because it \da{provides} only one effective coefficient at a macroscale point. \da{For accurate modeling, }the coarse-grid formulation needs multiple homogenized coefficients. \da{Multicontinuum approaches define different average states (continua) and introduce several effective properties per coarse-grid node. There are many works devoted to developing multicontinuum approaches for different applications \cite{rubinvstein1948question, barenblatt1960basic, arbogast1990derivation, panasenko2018multicontinuum}. Among them, one can mention a recently presented} multicontinuum homogenization method \cite{efendiev2023multicontinuum, chung2024multicontinuum, leung2024some}, \da{which allows us to rigorously derive multicontinuum models}. The main idea of \da{the} multicontinuum homogenization method is to derive a macroscopic model by formulating constraint cell problems and obtaining multicontinuum expansions. This method has already been \da{successfully applied to different problems} \cite{xie2024multicontinuum, ammosov2024multicontinuum, ammosov2023multicontinuum, efendiev2024multicontinuum}.

	
	\da{Note that one can also apply various} multiscale methods to numerically solve problems \da{with high-contrast properties}, such as the Multiscale Finite Volume Method (MsFVM) \cite{jenny2003multi, lunati2006multiscale, chaabi2024algorithmic}, the Multiscale Finite Element Method (MsFEM) \cite{hou1997multiscale, hou1999convergence, efendiev2009multiscale, jiang2017reduced}, the Generalized Multiscale Finite Element Method (GMsFEM) \cite{efendiev2013generalized, chung2014generalized, chung2014generalized1, chung2016adaptive, xie2025time} and the Constraint Energy Minimizing Generalized Multiscale Finite Element Method (CEM-GMsFEM) \cite{chung2018constraint, chung2021convergence, xie2024cem}. \da{In} these multiscale approaches\da{, one constructs} multiscale basis functions by \da{solving} local problems. \da{Then, the computed basis functions are used to obtain the coarse-grid problem.} \da{The multiscale methods demonstrated high efficiency in solving} time-fractional derivative problems with heterogeneous coefficient\da{s} \cite{tyrylgin2023computational, li2022partially, alikhanov2025multiscale}. \da{However, the obtained coarse-grid models are discrete and not in the form of macroscopic laws.}

	
	In \da{this} paper, we use the multicontinuum homogenization to the time-fractional diffusion-wave equation with high-contrast coefficient $\kappa$. \da{The main c}hallenges \da{lie in} formulating constraint cell problems, \da{deriving the multicontinuum model,} and \da{discretizing the} time-fractional derivative. We \da{formulate and solve} constraint cell problems that account for different averages and gradient effects. \da{Using the cell problems' solutions, }we obtain the multicontinuum expansions \da{and substitute them into} the variational \da{formulation} of the time-fractional diffusion-wave equation. Then\da{,} we \da{derive} a general multicontinuum model that includes macroscopic variables and homogenized coefficients. \da{For temporal discretization, we apply the fully discrete scheme \cite{sun2006fully} that can efficiently handle time-fractional derivatives.} To validate the efficiency of the proposed approach, we \da{consider two-dimensional model problems with high-contrast coefficients}. \da{We consider various orders of time-fractional derivatives, including the mixed case, where we have different orders in different continua.} \da{The} numerical results show that the \da{proposed multicontinuum approach provides high accuracy and can handle various heterogeneity types and time-fractional derivatives orders.}
	
	Our main contributions in this work are as follows. \da{First,} we \da{obtain} a macroscopic \da{model} for \da{a} time-fractional equation \da{with mixed time derivatives and zero-order spatial derivatives}. \da{Second,} we \da{derive} the multicontinuum \da{time-fractional diffusion-wave} model \da{using the multicontinuum homogenization approach}. \da{Finally,} we present numerical results to demonstrate the \da{efficiency} of the proposed \da{approach}.
	
	
	The rest of the paper is organized as follows. \da{In} the next section, we present preliminaries and derive \da{the} multicontinuum model for the  time-fractional zero-order equations. Section 3 is devoted to the derivation of \da{the} multicontinuum \da{time-fractional diffusion-wave} model. In \da{S}ection 4, we present numerical experiments to test the efficiency of our proposed \da{approach}. Finally, we present conclusions in Section 5.

	\section{Preliminaries}
	\da{In this section, we present preliminaries. First, we describe the fully discrete scheme used for approximating time-fractional derivatives. Then, we present the derivation of the multicontinuum equations for the case with mixed time-fractional derivatives and zero-order spatial terms using the multicontinuum homogenization method.}

	\subsection{Time-fractional diffusion-wave equation}
	\da{Let us} consider the following time-fractional diffusion-wave equation (TFDWE),
	\begin{equation}
		\frac{\partial^\alpha u}{\partial t^\alpha}-\nabla\cdot(\kappa\nabla u)=f,\quad \Omega\times[0,T],
		\label{TFDW}
	\end{equation}
	with the initial conditions
	\begin{equation}
		u(x,0)=\phi(x),\quad\frac{\partial u(x,0)}{\partial t}=\psi(x),
	\end{equation}
	and the boundary conditions
	\begin{equation}
		u(0,t)=\varphi_{1}(x),\quad u(L,t)=\varphi_{2}(x),\quad t>0,
	\end{equation}
	where $\kappa=\kappa(x)$ is a heterogeneous coefficient such that $\eta={\max\kappa/ \min\kappa} \gg 1$, $\phi(x),\psi(x),\varphi_{1}(x), \da{\text{ and }} \varphi_{2}(x)$ are known functions, $f(x,t)$ is a known source term, $T$ is the final time, $\Omega$ is a bounded domain in ${R}^{d}$ $(d=1,2,...)$ with boundary $\partial\Omega $, and 
	$$\frac{\partial^\alpha u}{\partial t^\alpha}=\frac{1}{\Gamma(2-\alpha)}\int\limits_0^t\frac{\partial^2u(x,s)}{\partial s^2}\frac{\mathrm{d}s}{(t-s)^{\alpha-1}},\quad1<\alpha<2$$
	is the Caputo fractional derivative of order $\alpha$ with $\Gamma$ denoting the gamma function. Note that if $\alpha = 1$, equation (\ref{TFDW}) represents a traditional diffusion equation. When $\alpha = 2$, equation (\ref{TFDW}) is considered as a traditional wave equation. For $1<\alpha<2$\da{,} equation (\ref{TFDW}) is known as the time-fractional diffusion-wave equation, \da{since} it can \da{describe intermediate processes} between \da{diffusion} ($\alpha = 1$) and wave ($\alpha = 2$) \da{phenomena}.
	
	To approximate the time-fractional derivative 	$\frac{\partial^\alpha u}{\partial t^\alpha}$, we first decompose the time interval $[0,T]$ evenly into $\da{N_t}$ time sub\da{intervals} $(t_{n},t_{n+1})$, $n=0,1,...,\da{N_t}-1$, with the time step size $\tau=\frac{T}{\da{N_t}}$. For ease of statement, we can suppose that $u(.,t_n)=\{u^n|0 \leq n\leq \da{N_t}\}$ is a grid function on $\Omega$, where $t_n=n\tau, \tau\textgreater0$. 
	\da{Let us introduce } the following notations
	$$\delta_{t}u^{n-\frac{1}{2}}=\frac{1}{\tau}(u^{n}-u^{n-1}),\quad u^{n-\frac{1}{2}}=\frac{1}{2}(u^{n}+u^{n-1}),$$
	where $u^{n-\frac{1}{2}}$ is the difference quotient of $u$ at the temporal points $t_n$ and $t_{n-1}$\da{,} and $u^{n-\frac{1}{2}}$ represents an average of $u$ on these two points.

	According to the fully discrete scheme in \cite{sun2006fully}, we can approximate the time-fractional derivative as follows
	\begin{equation}
		\frac{\partial^\alpha u}{\partial t^\alpha} \approx \frac{\tau^{1-\alpha}}{\Gamma(3-\alpha)}\biggl[a_0\delta_tu^{n-\frac{1}{2}}-\sum_{k=1}^{n-1}(a_{n-k-1}-a_{n-k})\delta_tu^{k-\frac{1}{2}}-a_{n-1}\psi\biggr]\\
		\label{time discrete}
	\end{equation}
	where $a_{k}=(k+1)^{2-\alpha}-k^{2-\alpha}(k\textgreater0).$ Sun and Wu \cite{sun2006fully} proved the solvability, unconditional stability\da{,} and $L_{\infty}$ convergence by the energy method. \da{Note that} the  convergence order for time is $O({\tau^{3-\alpha}})$. 
	
	For spatial approximation, we can consider the finite
	element approximation of (\ref{TFDW}). The weak formulation of (\ref{TFDW}) reads: \da{F}ind $u^{n}\in H_{0}^{1}(\Omega)$ such that for all $v\in H_{0}^{1}(\Omega)$
	\begin{equation}
		\begin{aligned}
			&\frac{\tau^{1-\alpha}}{\Gamma(3-\alpha)}\biggl[a_0 (\delta_tu^{n-\frac{1}{2}}, v) -\sum_{k=1}^{n-1}(a_{n-k-1}-a_{n-k})(\delta_tu^{k-\frac{1}{2}},v) -a_{n-1}(\psi,v) \biggr]\\
			&=-(\kappa\nabla u^{n-\frac{1}{2}}, \nabla v)+(f^{n-\frac{1}{2}},v), \quad 1\leq n\leq \da{N_t},
			\label{w}
		\end{aligned}
	\end{equation}	
	where $f^{n-\frac{1}{2}}=f(:,\frac{t_n+t_{n-1}}2).$

	\subsection{Multicontinuum time-fractional zero-order equations}

\da{Let us} consider the following time-fractional equation \da{with a zero-order spatial term}

\begin{equation}
	\frac{\partial^{\alpha_k} u}{\partial t^{\alpha_k}} \psi_{k}+A(x)u(x)=f(x),\quad  1<\alpha_k<2,
\label{zero}
\end{equation}	
where $A(x)$ is a scalar function with multiple scales and high contrast \da{symbol $k$ is the Einstein summation convention}. For example, let $A$ be a periodic  function comprises two different regions with highly diverse coefficients. We assume $\psi_k$ is the characteristic function for the region $\Omega_k$, called the $k$th continua $(k=1,...,N)$. \da{Here}, $\alpha_k$ are distinct \da{time derivative orders} in different regions. 

\begin{figure}[htbp]
	\centering
	\begin{tikzpicture}
		\fill[blue] (1.5,1) rectangle (2,1.5);
		\fill[blue] (4.6,1.1) rectangle (4.9,1.4);
		\fill[blue] (7.5,1) rectangle (8,1.5);
		\draw[step=0.5cm] (0,0) grid (3,3);
		\draw (4,0.5) rectangle (5.5,2);
		\draw (4.6,1.1) rectangle (4.9,1.4);
		\draw[step=0.5cm,shift={(6.5,0)}] (0,0) grid (2.5,2.5);
		\draw (10,0.5) rectangle (11.5,2);
		\draw[blue,opacity=0.5,->] (1.75,1) -- (3.95,0.5);
		\draw[blue,opacity=0.5,->] (1.75,1.5) -- (3.95,2);
		\draw[blue,opacity=0.5,->] (4.75,1.1) -- (6.45,0);
		\draw[blue,opacity=0.5,->] (4.75,1.4) -- (6.45,2.5);
		\draw[blue,opacity=0.5,->] (7.625,1) -- (9.9,0.5);
		\draw[blue,opacity=0.5,->] (7.625,1.5) -- (9.9,2);
		\node[anchor=north] at (1.5,-0.5) {$\Omega$};
		\node[anchor=north] at (4.5,-0.5) {$\omega$};
		\node[anchor=north] at (7.625,-0.5) {$R_{\omega}^+$};
		\node[anchor=north] at (10.75,-0.5) {$R_{\omega}$};
	\end{tikzpicture}
	\caption{Illustration of computation domain $\Omega$, coarse block $\omega$, RVE $R_{\omega}$, and the oversampling domain $R_{\omega}^+$.}
	\label{domain}
\end{figure}
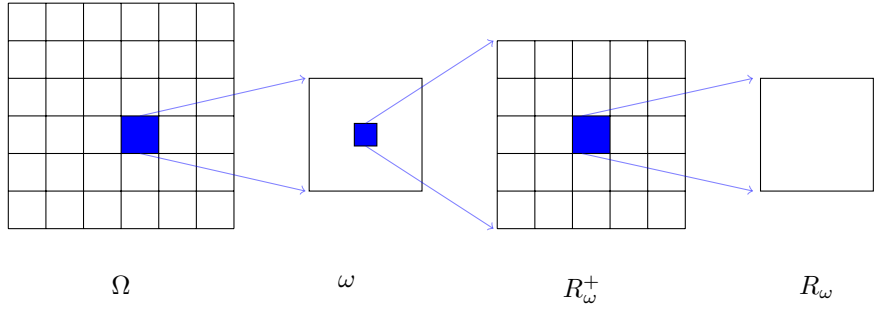

\da{Assume} that the problem is \da{defined} on a computational domain $\Omega$. \da{Moreover,} we divide the computational domain $\Omega$ into coarse blocks $\omega$'s. \da{Suppose that we can define} a representative volume (RVE) $R_{\omega}$ within \da{each} $\omega$, i.e., coarse\da{-}grid \da{elements} $\omega$ contain RVE\da{s}. Next, \da{we introduce an oversampled RVE} $R_{\omega}^+$, which is composed of several RVEs $R_{\omega}^p$\da{, where $R_{\omega}^{p_0} = R_{\omega}$}. We depict the relationships of $\Omega$, $\omega$, $R_{\omega}^+$\da{,} and $R_{\omega}$ in Figure \ref{domain}.
 
There are some following points to notice:
\begin{itemize}
	\item \da{Each} $R_{\omega}$ can \da{represent the whole} coarse-grid \da{element} $\omega$ \da{in terms of} heterogeneities;
	
	\item \da{Each} $R_{\omega}$ \da{contains} $N$ continua (components). \da{Moreover}, in each $R_{\omega}$, we \da{define} the characteristic function $\psi_j$ for continuum $j$, i.e., $\psi_{j}=\delta_{ij}$ within continuum $i$, where
	$$
	\delta_{ij}=\begin{cases} 1, \quad i=j,\\
		0, \quad i \neq j.
	\end{cases}
	$$
\end{itemize}

We define the expansion of $u$ in each RVE as
\begin{equation}
	u=\phi_iU_i,
	\label{U-e}
\end{equation}	
where $\phi_i$ is a microscopic function\da{,} and $U_i(x)$ is a smooth \da{macroscopic} function for each continuum $i$. Here, we use the Einstein summation convention. To \da{obtain} $\phi_i$, we introduce the following cell problems in each RVE  $R_{\omega}$. \da{Note that} we apply $y$ dependence to represent microscopic properties.

\begin{equation}
	\begin{aligned}
		A(y)\phi_{i}(y)&=D_{ij}\psi_{j} ,\quad \mathrm{in} ~R_{\omega},\\
		\int_{R_{\omega}}\phi_{i}\psi_{j}&=\delta_{ij}\int_{R_{\omega}}\psi_{j},\mathrm{~for~each~}j,
	\end{aligned}
\end{equation}
where $D_{ij}\psi_j=C_i \psi_i$ ($D_{ij}$ are constants). \da{One can see that} $A\phi_{i}=C_i\psi_i$ and
$$
C_j= \frac{\int_{R_{\omega}}\psi_{j}}{\int_{R_{\omega}}\psi_j^2A^{-1}}.
$$

Next, \da{let us derive} multicontinuum equations for $U_i$. We first write a weak \da{formulation of \eqref{zero}} for an \da{arbitrary} test function $v$
\begin{equation}
	\begin{aligned}
		\int_{\Omega}fv&=\int_{\Omega}\frac{\partial^{\alpha_k} u(x)}{\partial t^{\alpha_k}} \psi_{k}v(x)+\int_{\Omega} A(x)u(x)v(x)\\
		&=\sum_{\omega} \int_\omega\frac{\partial^{\alpha_k} u(x)}{\partial t^{\alpha_k}} \psi_{k}v(x)+\sum_{\omega} \int_\omega A(x)u(x)v(x)\\
		&\approx \sum_{\omega} \frac{|\omega|}{|R_{\omega}|} \int_{R_{\omega}}\frac{\partial^{\alpha_k} u(y)}{\partial t^{\alpha_k}} \psi_{k}v(y) + \sum_{\omega} \frac{|\omega|}{|R_{\omega}|}\int_{R_{\omega}}A(y)u(y)v(y).
	\end{aligned}
\label{weak}
\end{equation}	

\da{Next,} we substitute $u$ from (\ref{U-e}) into the equation (\ref{weak}), and \da{(}utilizing $v=\phi_i V_i$\da{)} we can approximate the first term of (\ref{weak}) as follows

\begin{equation}
	\begin{aligned}
		\int_{R_{\omega}}\frac{\partial^{\alpha_k} u(y)}{\partial t^{\alpha_k}} \psi_{k}v(y)& \approx \frac{\partial^{\alpha_k} U_i(x_\omega)}{\partial t^{\alpha_k}} V_j(x_\omega)\int_{R_{\omega}}\psi_k \phi_i(y) \phi_j(y)\\
		& = \frac{\partial^{\alpha_k} U_i(x_\omega)}{\partial t^{\alpha_k}} V_j(x_\omega) \int_{R_{\omega}} \psi_k \frac{C_i C_j}{A(y)A(y)}\psi_{i} \psi_{j},\\
	\end{aligned}
\end{equation}	
where $x_{\omega}$ is a mid-point of $R_{\omega}$. We will \da{omit} the microscopic dependence of macroscale variables (e.g.,$U_i$) and merely apply $U_i$ symbol
\begin{equation}
	\begin{aligned}
		\int_{R_{\omega}}\frac{\partial^{\alpha_k} u(y)}{\partial t^{\alpha_k}} \psi_{k}v(y)&=\frac{\partial^{\alpha_k} U_i}{\partial t^{\alpha_k}} V_j \int_{R_{\omega}} \psi_k \frac{C_i C_j}{A(y)A(y)}\psi_{i} \psi_{j}.
	\end{aligned}
\end{equation}	
\da{Let us denote}
\begin{equation}
\gamma_{ijk} = \int_{R_{\omega}} \psi_k \frac{C_i C_j}{A(y)A(y)}\psi_{i} \psi_{j}.
\end{equation}
\da{Note that $\gamma_{ijk}$ takes non-zero values only in the case of} $i=j=k$. \da{Therefore,} we can \da{represent it as follows} 
\begin{equation}
\gamma_i = C_i^2 \int_{R_{\omega}}\frac{\psi_{i}^3}{A^2}.
\end{equation}

The second term of (\ref{weak}) is approximated in a similar way
\begin{equation}
	\begin{aligned}
	\int_{R_{\omega}} A(y)u(y)v(x)&\approx U_i V_j\int_{R_{\omega}}A(y)\phi_i(y)\phi_j(y)\\
	&=U_i V_j \beta _{ij},
	\end{aligned}
\end{equation}	
where
\begin{equation}
	\begin{aligned}
    \beta_{ij}&= \int_{R_{\omega}} A(y) \phi_i(y) \phi_j(y) \\
    &=C_i \int_{R_{\omega}} \psi_i \phi_j \\
    &=C_i \delta_{ij} \int_{R_{\omega}} \psi_i.
	\end{aligned} 
	\label{alpha}
\end{equation}
\da{One can see that $\beta_{ij}$ takes non-zero values only in the case of} $i=j$. Thus, we \da{can represent it in the following form}	
\begin{equation}
\beta_{i}=	C_i \int_{R_{\omega}} \psi_i.
\end{equation}
From the above, we can obtain the following macroscopic equation\da{s} for (\ref{zero})
\begin{equation}
	\gamma_{ijk
	} \frac{\partial^{\alpha_{\da{k}}} U_i}{\partial t^{\alpha_{\da{k}}}} +\beta_{ij} U_i=b_j,\quad \da{j}=1,...,N,\quad 1<\alpha<2,	
\end{equation} 	
where $b_j = \int_{R_{\omega}} f \phi_j$.
	
Considering that $\gamma$ and $\alpha$ are diagonal matri\da{ces}, we have 

\begin{equation}
\gamma_i \frac{\partial^{\alpha_i} U_i}{\partial t^{\alpha_i}} +\beta_i U_i=b_i,\quad i=1,...,N,\quad 1<\alpha<2.	
\label{m-e}
\end{equation} 	

\da{One can see that the obtained macroscopic model is described by separate (not coupled) equations for each continuum}.


	\section{Multicontinuum homogenization}
	
	In this section, \da{we derive the multicontinuum time-fractional diffusion-wave model using the multicontinuum homogenization method}. We formulate constraint cell problems by \da{imposing} constraints for averages and gradient effects. These cell problems are \da{defined} in oversampled domains to reduce boundary effects. By applying the expansions to the variational formula, we can \da{rigorously derive} the multicontinuum model. \da{This section consists of two subsections. In the first subsection, we present the derivation for the regular time-fractional diffusion-wave equation. The second subsection considers the case with mixed time derivatives, where we have different time derivative orders in different continua.}

	
\subsection{Multicontinuum time-fractional diffusion-wave model}
   \da{Let us} consider the following variational formulation of the time-fractional diffusion-wave equation \da{\eqref{TFDW}}
	\begin{equation}
		\int_\Omega\frac{\partial^\alpha u}{\partial t^\alpha}v+\int_\Omega\kappa\nabla u\cdot\nabla v=\int_\Omega fv,\quad   1<\alpha<2, \quad \forall v\in H_0^1(\Omega).
	\label{var}
	\end{equation}
%

In the multicontinuum homogenization method, we \da{suppose} that \da{we can define} characteristic functions \da{of continua} in each RVE $R_\omega$ such that
$$U_i(x_\omega) \approx \frac{\int_{R_\omega}u\psi_{i}}{\int_{R_\omega}\psi_{i}}.$$
We assume that $U_i(x)$ \da{are} smooth function\da{s} denoting macroscopic variables, where $x_\omega$ is a \da{middle} point of $R_\omega$. Next, we present \da{the} derivation of the \da{multicontinuum} model for \da{the} time-fractional diffusion-wave equation.

\textbf{Expansion.} 

We expand the solution $u$ \da{over} macroscopic variables
 \begin{equation}
 u=\phi_iU_i+\phi_i^m\nabla_mU_i+\phi_i^{mn}\nabla_{mn}^2U_i+\cdots, 	
 \label{expansion}
 \end{equation}
 where $\nabla _m=\frac{\partial}{\partial x_m}$, $\nabla _{mn}^2=\frac{\partial^2}{\partial x_m \partial x_n}$. In this expansion, $U_i$ is the macroscopic quantity representing the average of the solution within $i$-th continuum\da{,} and $\phi_i,\phi_i^m,\phi_i^{mn},\ldots $ are auxiliary basis functions defined as the solutions of cell problems. For sake of simplicity, we will \da{neglect} the terms after the second term and only use the first two terms in (\ref{expansion}). \da{Therefore}, we can rewrite \da{the expansion of} $u$ as
 \begin{equation}
 	u\approx\phi_iU_i+\phi_i^m\nabla_mU_i,
 \end{equation}
where $U_i$ can be considered as a limit of $\int_{R_\omega}u\psi_i/\int_{R_\omega}\psi_i$.

\textbf{Cell problems.}	

Next, we introduce cell problems for auxiliary basis functions $\phi_i$ \da{and} $\phi_i^m$. These constrained cell problems are formulated by using Lagrange multipliers. To avoid boundary effects, we \da{define} \da{the} cell problems in oversampled regions $R_{\omega}^+$\da{, consisting} of $R_{\omega}^p$. We can minimize the local solution subject to different constraints \cite{chung2024multicontinuum, efendiev2023multicontinuum}.

The first cell problem considers different averages in each continuum and stands for the constants in the average behavior
 
 \begin{equation}
 	\begin{aligned}
 		\int_{R_\omega^+} \kappa\nabla\phi_i \cdot\nabla v
 		-\sum_{j,p}\frac{\beta_{ij}^p}{\int_{R_\omega^p} \psi_j^p}
 		\int_{R_\omega^p}\psi_j^p v &=0,\\
 		\int_{R_\omega^p}\phi_i\psi_j^p
 		&=\delta_{ij}\int_{R_\omega^p}\psi_j^p.
 	\end{aligned}
 \end{equation}
 The second cell problem considers gradient effects and imposes constraints to stand for the linear functions in the average behavior of each continuum
 
 \begin{equation}
 	\begin{aligned}
 		\int_{R_{\omega}^+} \kappa \nabla\phi_i^m \cdot\nabla v
 		-\sum_{j,p}\frac{\beta_{ij}^{mp}}{\int_{R_{\omega}^p}\psi_{j}^p}
 		\int_{R_{\omega}^p}\psi_{j}^p v &=0, \\
 		\int_{R_{\omega}^p} \phi_i^m \psi_j^p
 		&=\delta_{ij}\int_{R_{\omega}^p}(x_m-c_{m\da{j}})\psi_{j}^p,
 	\end{aligned}
 \end{equation}
 where $c_{m\da{j}}$ satisfy $\int_{R_{\omega}^{\da{p_0}}}(x_m-c_{m\da{j}})\psi_j^{\da{p_0}}=0.$
 
 \da{Note that we have} the following \da{estimates} for $\phi_i$ and $\phi_i^m$ \cite{efendiev2023multicontinuum}
 \begin{equation}
 \begin{aligned}
 ||\phi_i||&=O(1), \quad ||\nabla\phi_i||=O(\frac{1}{\epsilon}),\\
 ||\phi_i^m||&=O(\epsilon), \quad ||\nabla\phi_i^m||=O(1),
 \end{aligned}
 \label{scale}
 \end{equation}
 \da{where $\epsilon$ is a diameter of RVE.}
  
  \textbf{Integral localization.}	 
 
Utilizing the RVE concept, for any $v\in H_0^1$, we can get
 \begin{equation}
 	\begin{aligned}
 	\int_\Omega fv &= \int_\Omega\frac{\partial^\alpha u}{\partial t^\alpha}v+\int_\Omega\kappa\nabla u\cdot\nabla v\\
 	&=\sum_\omega\int_\omega\frac{\partial^\alpha u}{\partial t^\alpha}v+\sum_\omega\int_\omega\kappa\nabla u\cdot\nabla v\\
 	&\approx \sum_\omega \frac{|\omega|}{|R_\omega|}\int_{R_\omega}\frac{\partial^\alpha u}{\partial t^\alpha}v+\sum_\omega\frac{|\omega|}{|R_\omega|}\int_{R_\omega}\kappa\nabla u\cdot\nabla v,
 	\end{aligned}
 	\label{RVE}
 \end{equation}
where we assume that integrated average over RVE $R_\omega$ can \da{represent} the whole coarse block $\omega$ \da{in terms of heterogeneities}.

 \textbf{Substitution \da{into} the variational formulation.}	
 
 We can expand $u$ and $v$ as follows

 \begin{equation}
 	\begin{aligned}
 		u &\approx \phi_i U_i + \phi_i^m \nabla_m U_i, \\
 		v &\approx \phi_j V_j + \phi_j^n \nabla_n V_j,\\
 	\end{aligned}
 	\label{eq:uv_expansion}
 \end{equation}
 where $U_i$ and $V_j$ are smooth \da{macroscopic} functions.

Substituting \eqref{eq:uv_expansion} into \eqref{RVE}, we can obtain the following \da{approximation of the variational formulation}
\begin{equation}
 	\begin{aligned}
 		&\sum_{R_\omega}\frac{|\omega|}{|R_\omega|}
 		\Big(\int_{R_\omega}\frac{\partial^{\alpha} (\phi_i U_i + \phi_i^m \nabla_m U_i)}{\partial t ^{\alpha}} (\phi_j V_j + \phi_j^n \nabla_n V_j)+\\ &\sum_{R_\omega}\frac{|\omega|}{|R_\omega|}\int_{R_\omega} \kappa \nabla (\phi_i U_i + \phi_i^m \nabla_m U_i)\cdot\nabla (\phi_j V_j + \phi_j^n \nabla_n V_j)=\\
        &\sum_{R_\omega}\frac{|\omega|}{|R_\omega|}\int_{R_\omega}\{(f ,\phi_j V_j) + (f, \phi_j^n \nabla_n V_j)\}.
 	\end{aligned}
 	\label{va}
\end{equation}

\textbf{Macroscopic model.}	

Due to the smoothness of the macroscopic variables, the variations of $U_i$ and $\nabla_{m}U_i$ are small compared to the variations of $\phi_{i}$ and $\phi_{i}^m$. \da{Therefore,} we assume $\int_{R_{\omega}} \kappa \nabla (\phi_{i} U_i) \cdot \nabla v \approx \int_{R_{\omega}} \kappa \nabla (\phi_{i}) U_i \cdot \nabla v$ and 
$\int_{R_{\omega}} \kappa \nabla (\phi_i^m \nabla_m U_i) \cdot \nabla v \approx \int_{R_{\omega}} \kappa \nabla (\phi_i^m) \nabla_m U_i \cdot \nabla v$. As before (in \da{the} zero-order \da{case}), we take the macroscopic variables out of the integrals over $R_\omega$. 

We can approximate the first term on the left-hand side of (\ref{va}) \da{as follows}
\begin{equation}
	\begin{aligned}
		&\sum_{R_\omega}\frac{|\omega|}{|R_\omega|}
		\Big(\int_{R_\omega}\frac{\partial^{\alpha} (\phi_i U_i + \phi_i^m \nabla_m U_i)}{\partial t ^{\alpha}} (\phi_j V_j + \phi_j^n \nabla_n V_j)\Big)=\\
		&\sum_{\omega}\frac{|\omega|}{|R_{\omega}|}(\int_{R_{\omega}}\phi_{i}\phi_{j})\frac{\partial^{\alpha} }{\partial t ^{\alpha}}U_{i}V_{j}+\sum_{\omega}\frac{|\omega|}{|R_{\omega}|}(\int_{R_{\omega}}\phi_{i}\phi_{j}^{n})\frac{\partial^{\alpha} }{\partial t ^{\alpha}}U_{i}\nabla_{n}V_{j}+ \\
		&\sum_{\omega}\frac{|\omega|}{|R_{\omega}|}(\int_{R_{\omega}}\phi_{i}^{m}\phi_{j})\frac{\partial^{\alpha} }{\partial t ^{\alpha}}\nabla_{m}U_{i}V_{j}+\sum_{\omega}\frac{|\omega|}{|R_{\omega}|}(\int_{R_{\omega}}\phi_{i}^{m}\phi_{j}^{n})\frac{\partial^{\alpha} }{\partial t ^{\alpha}}\nabla_{m}U_{i}\nabla_{n}V_{j}\\
		&\approx\sum_{\omega}\frac{|\omega|}{|R_{\omega}|}(\int_{R_{\omega}}\phi_{i}\phi_{j})\frac{\partial^{\alpha} }{\partial t ^{\alpha}}U_{i}V_{j}.
	\end{aligned}
	\label{f_l}
\end{equation}
\da{In addition, let us introduce the following notations corresponding to the coefficients of the terms in \eqref{f_l}}
$$
\begin{aligned}
	C_{ji} &= \int_{R_{\omega}}\phi_{i}\phi_{j}, \quad 
	C_{ji}^{n} = \int_{R_{\omega}}\phi_{i}\phi_{j}^{n}, \\
	C_{ji}^{m} &= \int_{R_{\omega}}\phi_{i}^{m}\phi_{j}, \quad 
    C_{ji}^{mn} = \int_{R_{\omega}}\phi_{i}^{m}\phi_{j}^{n}.
\end{aligned}
$$

\da{Note that} in the last step of (\ref{f_l}), we cancel the second, third\da{,} and fourth terms. The reason for these is that $\phi_{j}^{n}$ is of \da{the} order $O(\epsilon)$, while $\phi_{j}$'s order is \da{of the order} $O(1)$.

The second term on the left-hand side of (\ref{va}) can be written as 
\begin{equation}
	\begin{aligned}	
&\sum_{R_\omega}\frac{|\omega|}{|R_\omega|}\int_{R_\omega} \kappa \nabla (\phi_i U_i + \phi_i^m \nabla_m U_i)\cdot\nabla (\phi_j V_j + \phi_j^n \nabla_n V_j)=\\
& \sum_{R_\omega}\frac{|\omega|}{|R_\omega|}\left\{ ({\int_{R_{\omega}} \kappa  \nabla \phi_i \cdot \nabla \phi_j})U_{i}V_{j} + ({\int_{R_{\omega}} \kappa  \nabla \phi_i \cdot \nabla \phi_j^n})U_{i}\nabla_{n}V_{j} \right\}+ \\
& \sum_{R_\omega}\frac{|\omega|}{|R_\omega|}\left\{({\int_{R_{\omega}} \kappa  \nabla \phi_i^m \cdot \nabla \phi_j})\nabla_{m}U_{i}V_{j} + ({\int_{R_{\omega}} \kappa  \nabla \phi_i^m \cdot \nabla \phi_j^n})\nabla_{m}U_{i}\nabla_{n}V_{j}\right\}.
	\end{aligned}
\end{equation}
\da{Accordingly, we introduce the following notations}
$$
\begin{aligned}
	B_{ji}& = {\int_{R_{\omega}} \kappa  \nabla \phi_i \cdot \nabla \phi_j, \quad 
	B_{ji}^n = {\int_{R_{\omega}} \kappa  \nabla \phi_i \cdot \nabla \phi_j^n}}, \\
	B_{ji}^{m} &= {\int_{R_{\omega}} \kappa  \nabla \phi_i^m \cdot \nabla \phi_j}, \quad 
	B_{ji}^{mn} = {\int_{R_{\omega}} \kappa  \nabla \phi_i^m \cdot \nabla \phi_j^n}.
\end{aligned}
$$

Applying the formula (\ref{scale}), we can \da{obtain the following} estimate\da{s for} $B_{ji}$, $C_{ji}$, $B_{ji}^m$, \da{and} $B_{ji}^{mn}$
$$
B_{ji}=O(\frac{|R_\omega|}{\epsilon^2}), \quad B_{ji}^m=O(\frac{|R_\omega|}{\epsilon}), \quad B_{ji}^{mn}=O(|R_\omega|), \quad C_{ji}=O(|R_\omega|).
$$

Then\da{,} we can define the scaled effective properties as follow\da{s}
\begin{equation}
	\widehat{B_{ji}}=\frac{\epsilon^2}{|R_\omega|}B_{ji},\quad \widehat{B_{ji}^m}=\frac{\epsilon}{|R_\omega|}B_{ji}^m,\quad \widehat{B_{ji}^{mn}}=\frac{1}{|R_\omega|}B_{ji}^{mn},\quad \widehat{C_{ji}}=\frac{1}{|R_\omega|}C_{ji}.
\end{equation}

Hence, with these scalings, we can get
\begin{equation}
	\begin{aligned}
		&\int_{\Omega}\frac{\partial^{\alpha}u}{\partial t^{\alpha}}v+\int_{\Omega}\kappa\nabla u\cdot\nabla v\approx\int_{\Omega}\widehat{C_{ji}}\frac{\partial^{\alpha}U_{i}}{\partial t^{\alpha}}V_{j}+\int_{\Omega}\widehat{B_{ji}^{mn}}\nabla_{m}U_{i}\nabla_{n}V_{j}+ \\
		&\frac{1}{\epsilon}\int_{\Omega}\widehat{B_{ji}^{m}}\nabla_{m}U_{i}V_{j}+\frac{1}{\epsilon}\int_{\Omega}\widehat{B_{ji}^{m}}U_{i}\nabla_{m}V_{j}+\frac{1}{\epsilon^{2}}\int_{\Omega}\widehat{B_{ji}}U_{i}V_{j}.
	\end{aligned}
	\label{w_s}
\end{equation}
By utilizing the integration by parts, \da{one can show that} the sum of the third and fourth term\da{s} on the right-hand side of (\ref{w_s}) is negligible \cite{efendiev2023multicontinuum}.

Thus, we \da{obtain} the following multicontinuum time-fractional diffusion-wave model (in strong form)

\begin{equation}
	\widehat{C_{ji}}\frac{\partial^\alpha U_i}{\partial t^\alpha}-\nabla_n(\widehat{B_{ji}^{mn}}\nabla_mU_i)+\frac{1}{\epsilon^2}\widehat{B_{ji}}U_i=f_j.
	\label{mc}
\end{equation}
where $f_j=\da{\frac{1}{|R_{\omega}|}} \int_{R_\omega}f \phi_j$. \da{One can see that the reaction terms are dominant unless we have large diffusions, i.e., high contrast.}

\subsection{Multicontinuum time-fractional diffusion-wave model with mixed time derivatives}

\da{Let us now} consider the following time-fractional diffusion-wave equation \da{with mixed time derivatives}

\begin{equation}
	\frac{\partial^{\alpha_p} u}{\partial t^{\alpha_p}}\psi_{p}-\nabla\cdot(\kappa\nabla u)=f,\quad 1<\alpha_p<2,
	\quad \Omega\times[0,T],
	\label{MTFDW}
\end{equation}
where $\kappa$ is \da{a high-contrast} heterogeneous coefficient, and $\psi_{p}$ is the characteristic function for the \da{subdomain} $\Omega_p$\da{,} i.e.\da{,} the $p$th continuum \da{(}$p=1,...,N$\da{)}. \da{Here, we suppose the Einstein summation convention over $p$ indices.}

Applying the RVE, we first obtain a weak form of (\ref{MTFDW}) for an \da{arbitrary} test function $v \in H^1_0$ as follows
 \begin{equation}
	\begin{aligned}
		\int_\Omega fv &= \int_\Omega\frac{\partial^{\alpha_p} u}{\partial t^{\alpha_p}}\psi_{p}v+\int_\Omega\kappa\nabla u\cdot\nabla v\\
		&=\sum_\omega\int_\omega\frac{\partial^{\alpha_p} u}{\partial t^{\alpha_p}}\psi_{p}v+\sum_\omega\int_\omega\kappa\nabla u\cdot\nabla v\\
		&\approx \sum_\omega \frac{|\omega|}{|R_\omega|}\int_{R_\omega}\frac{\partial^{\alpha_p} u}{\partial t^{\alpha_p}}\psi_{p}v+\sum_\omega\frac{|\omega|}{|R_\omega|}\int_{R_\omega}\kappa\nabla u\cdot\nabla v,
	\end{aligned}
	\label{weak1}
\end{equation}

The cell problems are \da{the} same as \da{for} the regular time-fractional diffusion-wave equation. \da{Thus,} we can define the following multicontinuum expansion of $u$ and $v$

\begin{equation}
	\begin{aligned}
		u &\approx \phi_i U_i + \phi_i^m \nabla_m U_i, \\
		v &\approx \phi_j V_j + \phi_j^n \nabla_n V_j.\\
	\end{aligned}
	\label{eq:uv_expansion1}
\end{equation}
Then\da{,} we can \da{approximate} (\ref{weak1}) as \da{follows}

\begin{equation}
	\begin{aligned}
		&\sum_{R_\omega}\frac{|\omega|}{|R_\omega|}
		\Big(\int_{R_\omega}\frac{\partial^{{\alpha_p}} (\phi_i U_i + \phi_i^m \nabla_m U_i)}{\partial t ^{{\alpha_p}}} \psi_{p}(\phi_j V_j + \phi_j^n \nabla_n V_j)+\\ &\sum_{R_\omega}\frac{|\omega|}{|R_\omega|}\int_{R_\omega}\nabla (\phi_i U_i + \phi_i^m \nabla_m U_i)\cdot\nabla (\phi_j V_j + \phi_j^n \nabla_n V_j)=\\
		&\sum_{R_\omega}\frac{|\omega|}{|R_\omega|}\{(f ,\phi_j V_j) + (f, \phi_j^n \nabla_n V_j)\}.
	\end{aligned}
	\label{va1}
\end{equation}

The first term on the left-hand side of (\ref{va1}) can be rewritten \da{in the following form} (for brevity, we will \da{omit} $\sum_{\omega}\frac{|\omega|}{|R_{\omega}|}$)

\begin{equation}
	\begin{aligned}
		&\int_{R_\omega}\frac{\partial^{\alpha_p} (\phi_i U_i + \phi_i^m \nabla_m U_i)}{\partial t ^{\alpha_p}} \psi_{p}(\phi_j V_j + \phi_j^n \nabla_n V_j)\\
		&=C_{jip}\frac{\partial^{\alpha_p} }{\partial t ^{\alpha_p}}U_{i}V_{j}+C_{jip}^{n}\frac{\partial^{\alpha_p} }{\partial t ^{\alpha_p}}U_{i}\nabla_{n}V_{j}+ 
		C_{jip}^{m}\frac{\partial^{\alpha_p} }{\partial t ^{\alpha_p}}\nabla_{m}U_{i}V_{j}+	C_{jip}^{mn}\frac{\partial^{\alpha_p} }{\partial t ^{\alpha_p}}\nabla_{m}U_{i}\nabla_{n}V_{j}\\
		&\approx C_{jip}\frac{\partial^{\alpha_p} }{\partial t ^{\alpha_p}}U_{i}V_{j}.
	\end{aligned}
	\label{f_l1}
\end{equation}

where 
$$
\begin{aligned}
	C_{jip} &= \int_{R_{\omega}}\phi_{i}\phi_{j}\psi_{p}, \quad 
	C_{jip}^{n} = \int_{R_{\omega}}\phi_{i}\phi_{j}^{n}\psi_{p}, \\
	C_{jip}^{m} &= \int_{R_{\omega}}\phi_{i}^{m}\phi_{j}\psi_{p}, \quad 
	C_{jip}^{mn} = \int_{R_{\omega}}\phi_{i}^{m}\phi_{j}^{n}\psi_{p}.
\end{aligned}
$$

The second term on the left-hand side of (\ref{va1}) is \da{the} same as \da{in} the regular case. \da{Therefore,} we can write it as follows

\begin{equation}
		\int_{R_\omega} \kappa\nabla u\cdot\nabla v\approx
		\int_{R_\omega}	B_{ji}U_{i}V_{j} + \int_{R_\omega} B_{ji}^nU_{i}\nabla_{n}V_{j} + 
		 \int_{R_\omega} 	B_{ji}^{m}\nabla_{m}U_{i}V_{j} + \int_{R_\omega}	B_{ji}^{mn}\nabla_{m}U_{i}\nabla_{n}V_{j},
\end{equation}
where 
$$
\begin{aligned}
	B_{ji}& = {\int_{R_{\omega}} \kappa  \nabla \phi_i \cdot \nabla \phi_j, \quad 
		B_{ji}^n = {\int_{R_{\omega}} \kappa  \nabla \phi_i \cdot \nabla \phi_j^n}}, \\
	B_{ji}^{m} &= {\int_{R_{\omega}} \kappa  \nabla \phi_i^m \cdot \nabla \phi_j}, \quad 
	B_{ji}^{mn} = {\int_{R_{\omega}} \kappa  \nabla \phi_i^m \cdot \nabla \phi_j^n}.
\end{aligned}
$$

\da{Note that} we \da{can obtain} the \da{estimate} $C_{jip}=O(|R_\omega|)$. Then\da{,} we \da{have the following scaled effective property} $\widehat{C_{jip}}=\frac{1}{|R_\omega|}C_{jip}$. Eventually, we \da{obtain} the following \da{multicontinuum model} (in strong form)
\begin{equation}
	\widehat{C_{jip}}\frac{\partial^{\alpha_{\da{p}}} U_i}{\partial t^{\alpha_{\da{p}}}}-\nabla_n(\widehat{B_{ji}^{mn}}\nabla_mU_i)+\frac{1}{\epsilon^2}\widehat{B_{ji}}U_i=f_j.
	\label{mc1}
\end{equation}
\da{Thus, each equation contains a linear combination of time-fractional derivatives of all orders.}

	\section{Numerical results}

	
In this section, we \da{conduct numerical experiments to check the efficiency of the proposed multicontinuum homogenization approach. For this purpose, we consider two-dimensional model problems in different high-contrast heterogeneous media. We consider various time-fractional derivative orders, including the mixed derivatives case.} 

\da{As a computational domain, we set} $\Omega = \Omega_1 \cup \Omega_2 = [0,1] \times [0,1]$, \da{where $\Omega_1$ is a low conductive region, and $\Omega_2$ is a highly conductive region.} 
\da{To construct a fine grid,} we divide $\Omega$ into $400 \times 400$ square \da{blocks}, 
We partition $\Omega$ evenly into $M \times M$ square \da{blocks (each divided into two triangles) to construct a coarse grid}. \da{Therefore,} the coarse \da{grid} size is $H = 1 / M$. For sake of simplicity, each whole coarse-grid \da{block} will be considered as an RVE \da{itself}. The oversampled RVE $R^{+}_{\omega}$ is \da{constructed} as an extension of \da{the target block} by $l$ layers of coarse-grid \da{blocks}. \da{To estimate the accuracy of the proposed approach,} we \da{use} the \da{following} relative $L^2$-error \da{norms}
\begin{equation}
	e^{(i)}(t)=\sqrt{\frac{\sum_K|\frac{1}{|K|}\int_KU_i(x,t)dx-\frac{1}{|K\cap\Omega_i|}\int_{K\cap\Omega_i}u(x,t)dx|^2}{\sum_K|\frac{1}{|K\cap\Omega_i|}\int_{K\cap\Omega_i}u(x,t)dx|^2}},
\end{equation}
where $i=1,2$, and $K$ denotes the RVE, which is taken to be $\omega$.  

\da{We consider two coarse grid sizes:} $1/20$ \da{and} $1/40$. \da{We determine} the number of oversampling layers $l$ using the formula $l=\lceil-2\log(H)\rceil$ \da{\cite{efendiev2023multicontinuum}}. Hence, we \da{have} $l=6$ for $H=1/20$ and $l=8$ for $H=1/40$. We utilize the finite element method \da{with standard piece-wise linear basis functions for spatial approximation of both fine-grid and coarse-grid models}. For temporal approximation, we apply proposed the fully discrete scheme (\ref{time discrete}). We take the time step $\tau = 0.02$ and \da{set} the final time $T=1$. 

\da{The numerical implementation is based on the open-source computational package FEniCS \cite{logg2012automated}. The visualization of the obtained numerical results was performed using the ParaView software \cite{ahrens2005paraview}.}

\subsection{Example 1: Crossed field}

\begin{figure}[hbt!]
	\centering
	\includegraphics[width=0.29\textwidth]{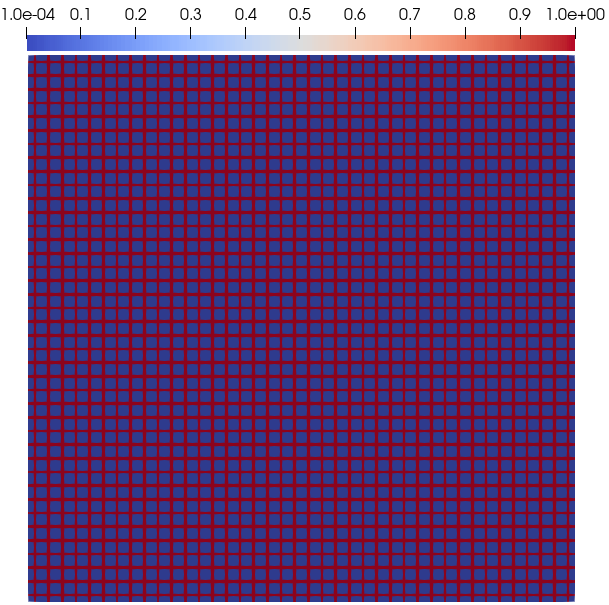}
	\caption{The coefficient $\kappa$ ($\Omega_1$: blue regions; $\Omega_2$: red regions). \da{Crossed field.}}
	\label{kappa}
\end{figure}

\da{First, let us consider the heterogeneous coefficient $\kappa$ with a crossed structure (see Figure \ref{kappa}). Moreover, we suppose that $\kappa$ posesses a high contrast and can be represented in the following form}
\begin{equation}
	\kappa(x)=
	\begin{cases}
		\begin{aligned}
			10^{-4}, \quad &x \in \Omega_1,\\
			1, \quad &x \in \Omega_2.
		\end{aligned}
	\end{cases}
	\label{eq:kappa}
\end{equation} 

The source term $f$ is given as $f(x)=e^{-40((x_1 -0.5)^2+(x_2 -0.5)^2)}$ for any $x=(x_1,x_2)\in \Omega$. We set the initial condition $u_0 = u(x,0)=0$\da{,} and the following Dirichlet boundary condition
\begin{equation}
	u(x,t) = 0, \quad x \in \partial \Omega.
\end{equation}

\da{We consider two cases of time-fractional derivative orders. In the first case, we set the same orders for both continua (regions). The second case corresponds to the mixed time derivatives, where we have different orders in different continua.}

\subsubsection{Case 1: \da{R}egular \da{time derivatives}}

\begin{figure}[hbt!]
	\centering
	\includegraphics[width=0.29\textwidth]{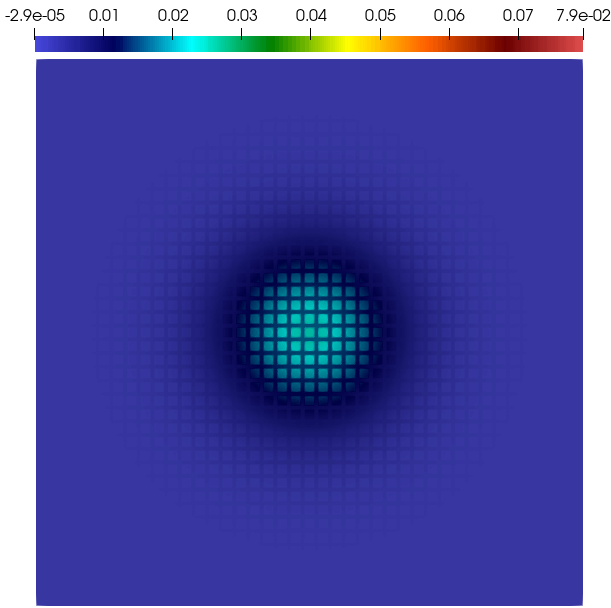}
	\includegraphics[width=0.29\textwidth]{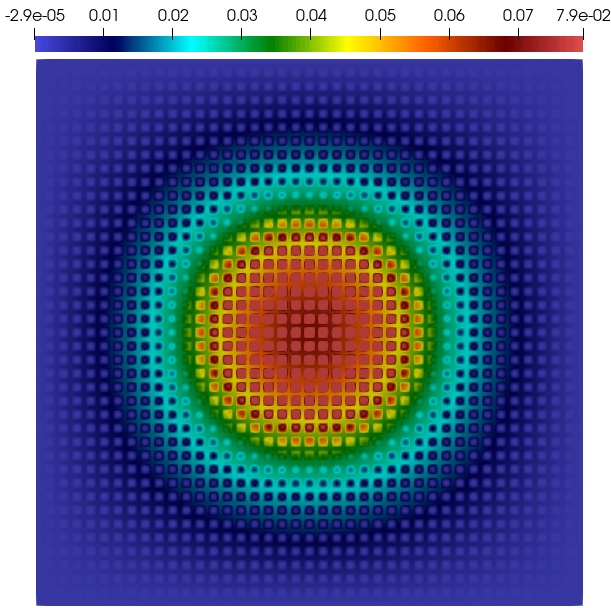}
	\includegraphics[width=0.29\textwidth]{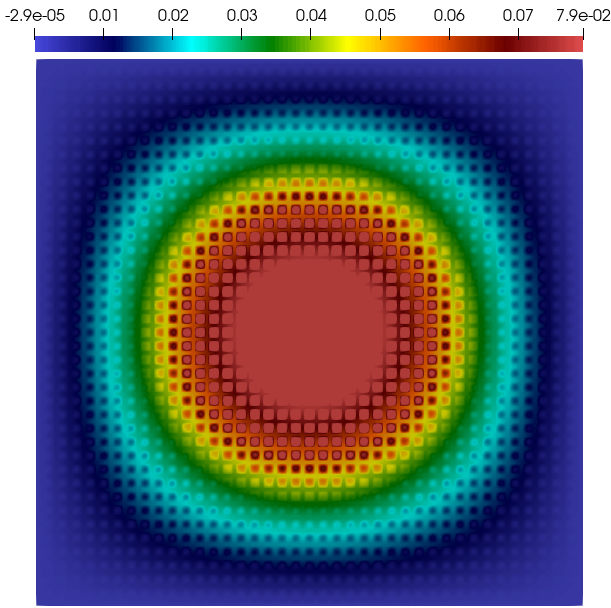}
	\caption{The fine-grid reference solution with $\alpha = 1.5$ at $t=0.1,0.5,1$ for Case 1 in Example 1 (from left to right).}
	\label{1.5_f}
\end{figure}

\da{Let us consider the case with the same time-fractional derivatives in both continua. Figure \ref{1.5_f} depicts distributions of the fine-grid solution for $\alpha = 1.5$ at different time steps. One can see how the solution field grows from the middle of the domain due to the source term. We can observe the high-contrast diffusion coefficient's influence on the distribution of the solution. The solution field propagates faster in the high-conductive channels while it moves slowly in the low-conductive subregions. In general, the obtained solution corresponds to the simulated process.}


\begin{figure}[hbt!]
	\centering
	\includegraphics[width=0.29\textwidth]{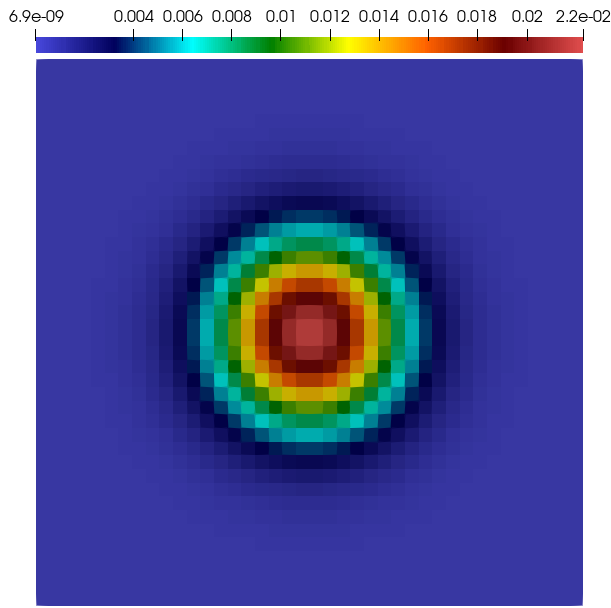}
	\includegraphics[width=0.29\textwidth]{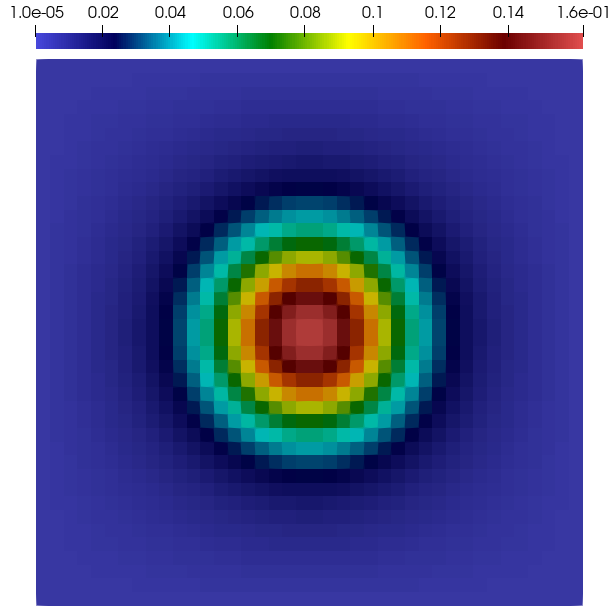}
	\includegraphics[width=0.29\textwidth]{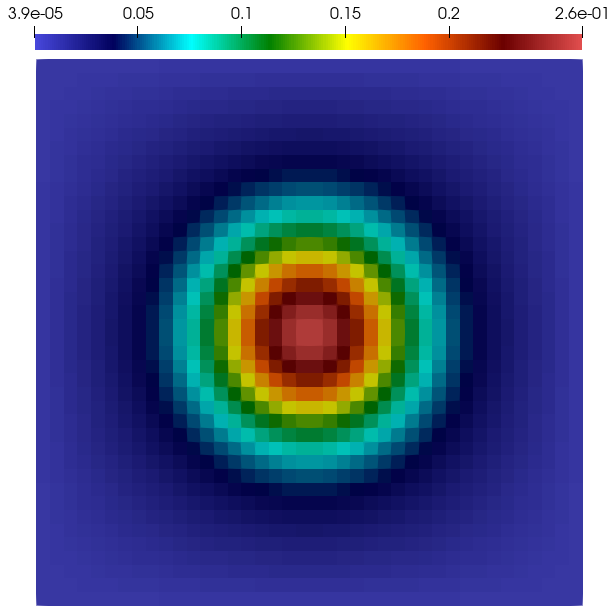}\\
	\includegraphics[width=0.29\textwidth]{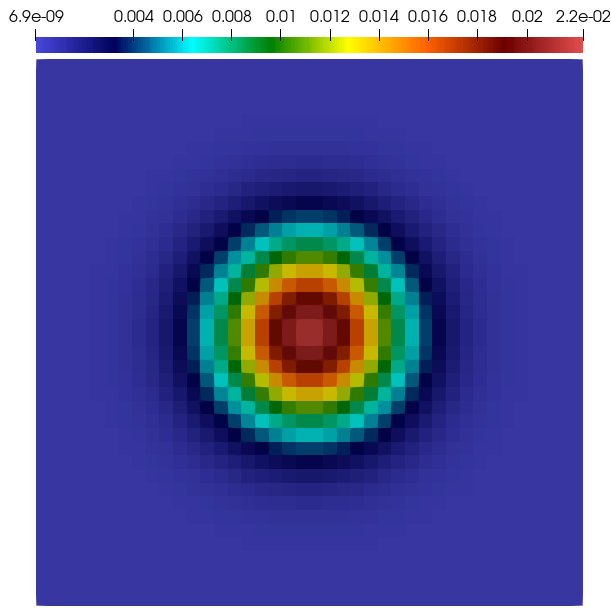}
	\includegraphics[width=0.29\textwidth]{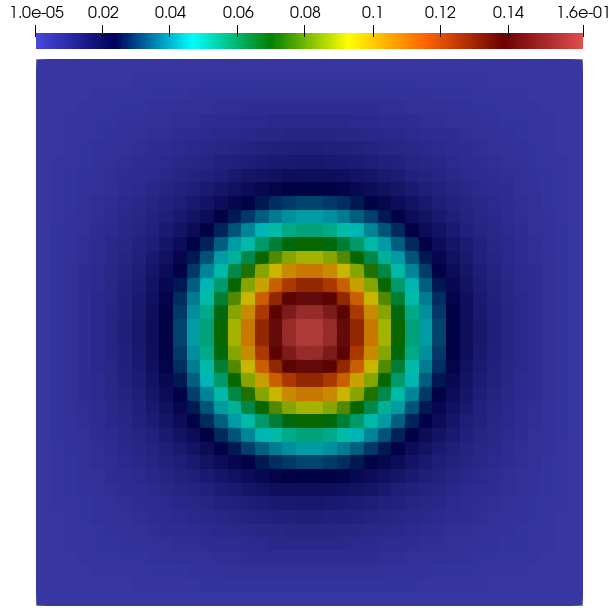}
	\includegraphics[width=0.29\textwidth]{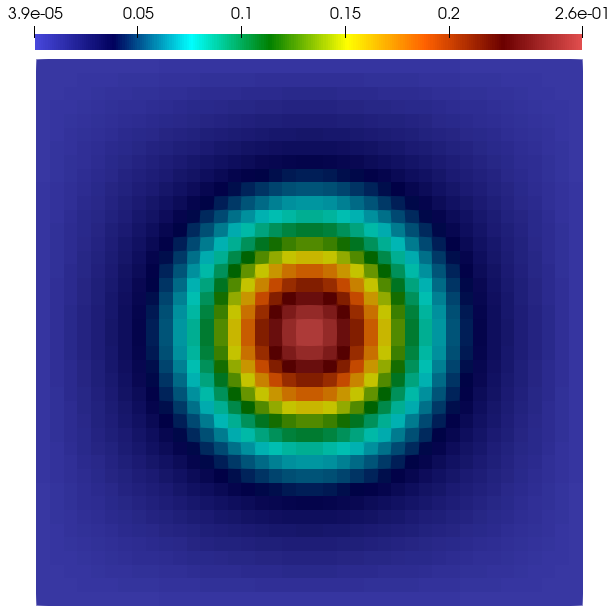}
	\caption{Average solution with $\alpha=1.5$ at $t=0.1,0.5,1$ for Case 1 in Example 1 (from left to right)\da{, $H = 1/40$}. First row: reference averaged solution in $\Omega_1$. Second row: multiscale solution in $\Omega_1$.}
	\label{1.5_1}
\end{figure}

\begin{figure}[hbt!]
	\centering
	\includegraphics[width=0.29\textwidth]{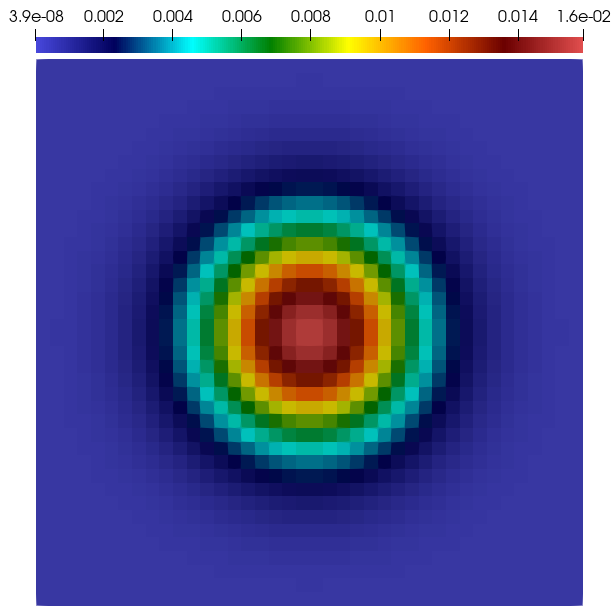}
	\includegraphics[width=0.29\textwidth]{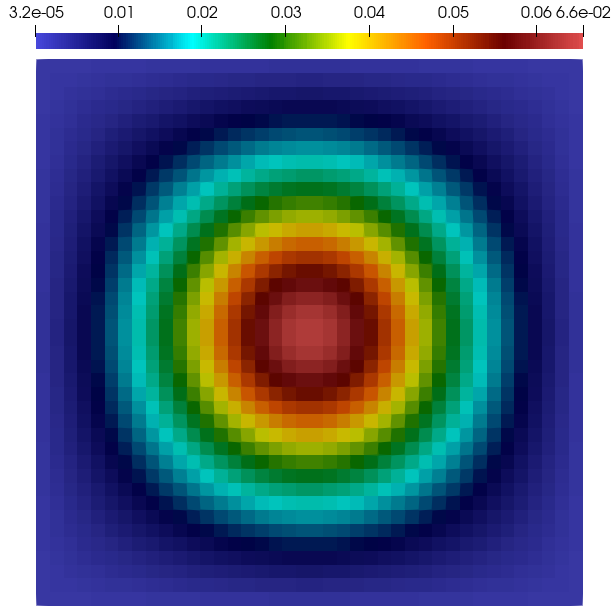}
	\includegraphics[width=0.29\textwidth]{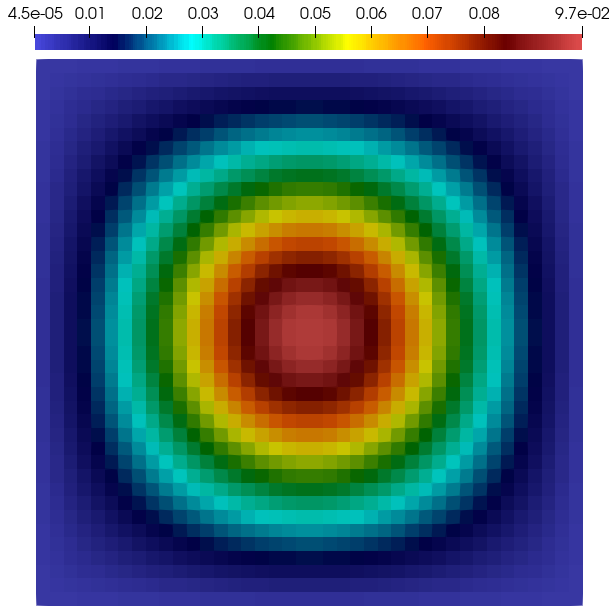}\\
	\includegraphics[width=0.29\textwidth]{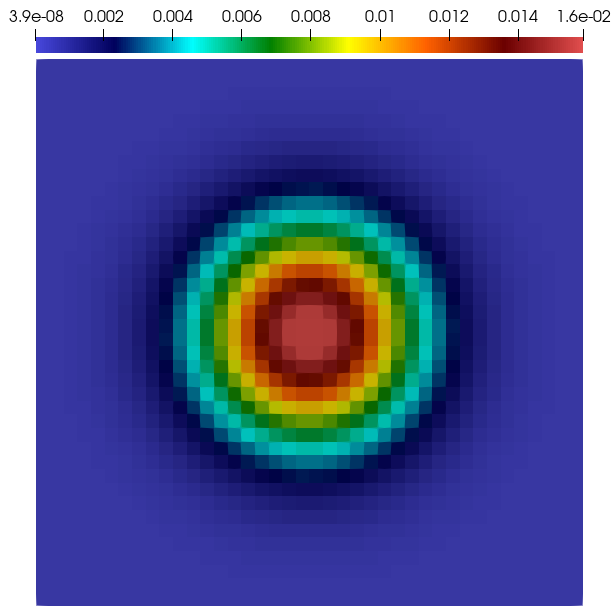}
	\includegraphics[width=0.29\textwidth]{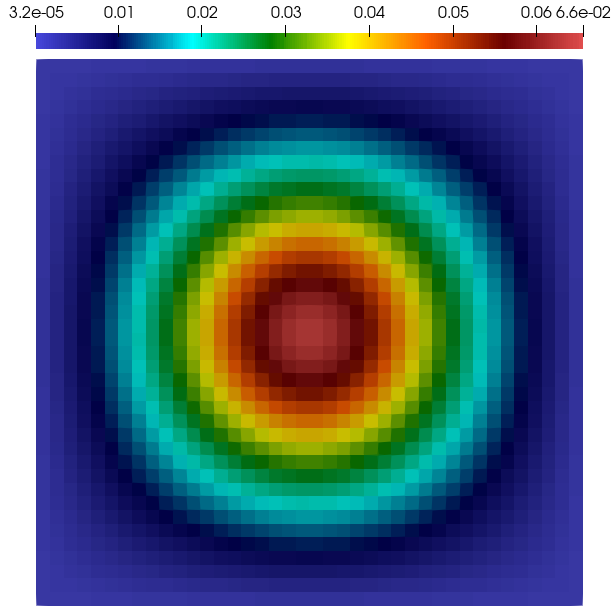}
	\includegraphics[width=0.29\textwidth]{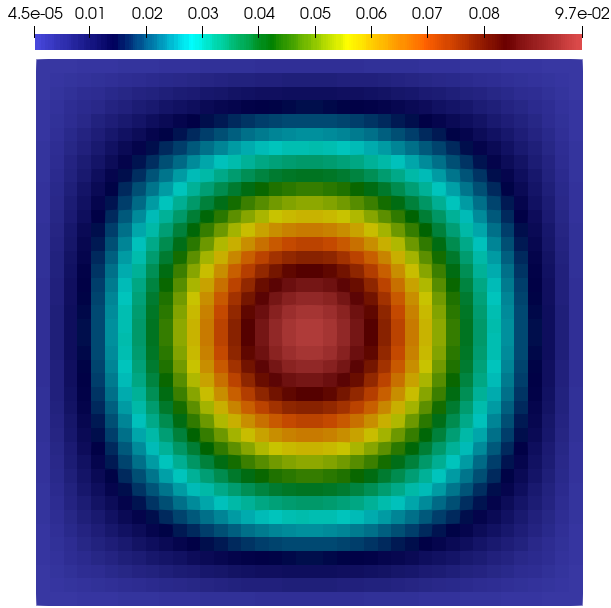}
	\caption{Average solution with $\alpha=1.5$ at $t=0.1,0.5,1$ for Case 1 in Example 1 (from left to right)\da{, $H = 1/40$}. First row: reference averaged solution in $\Omega_2$. Second row: multiscale solution in $\Omega_2$.}
	\label{1.5_2}
\end{figure}

\da{Figures \ref{1.5_1}-\ref{1.5_2} presents distributions of the average solutions with $H = 1 / 40$ in $\Omega_1$ and $\Omega_2$ subregions, respectively. From top to bottom, we depict the reference and multiscale average solutions. One can see that the solutions are very similar, which indicates that our proposed multicontinuum approach can approximate the reference solution with high accuracy. In terms of the simulated process, we observe that the average solution field in $\Omega_2$ ($U_1$) diffuses much faster than those in $\Omega_2$ ($U_2$), as expected. Moreover, one can see that the propagation process is isotropic.}

%
%
%




\da{Next, let us consider the relative errors of the multiscale solutions. Tables \ref{1.1error}, \ref{1.5error}, and \ref{1.9error} present the relative $L^2$ errors at different time steps for $\alpha = 1.1$, $\alpha = 1.5$, and $\alpha = 1.9$, respectively. In all the tables, we present the errors for $H = 1/20$ on the left and the errors for $H = 1/40$ on the right. One can see that the errors of the multiscale solutions for all $\alpha$ cases are minor. Moreover, we observe the error reduction with a decrease in the coarse grid size.}

\begin{table}[hbt!]
\caption{Relative errors at \da{different time steps} with $\alpha=1.1$ for Case 1 in Example 1. Left: $H=1/20$ and $l=6.$ Right: $H=1/40$ and $l=8.$}		
\begin{minipage}[c]{0.5\textwidth}
	\centering
	\begin{tabular}{c c c}
		\hline
	$t$ & $e^{(1)}(t)$ & $e^{(2)}(t)$ \\
		\hline
		0.1 & 4.3610\% & 1.6115\% \\
		
		0.2 & 3.0644\% & 1.8407\% \\
	
		0.3 & 2.5377\% & 1.6898\% \\
		
		0.4 & 2.4026\% &  1.5983\% \\
	
		0.5 & 2.4039\% & 1.5183\% \\
	
		0.6 & 2.4435\% & 1.4562\% \\
		
		0.7 & 2.4864\% & 1.4169\% \\
		
		0.8 & 2.5212\% & 1.3965\% \\
		
		0.9 & 2.5457\% & 1.3891\%\\
		
		1.0 & 2.5610\% & 1.3895\% \\
		\hline
	\end{tabular}
\end{minipage}
\begin{minipage}[c]{0.5\textwidth}
	\centering
	\begin{tabular}{c c c}
		\hline
		$t$ & $e^{(1)}(t)$ & $e^{(2)}(t)$ \\
		\hline
		0.1 & 2.6539\% & 0.9813\% \\
		
		0.2 & 1.3138\% & 0.7866\% \\
		
		0.3 & 0.7241\% & 0.6480\% \\
		
		0.4 & 0.5579\% & 0.5979\% \\
		
		0.5 & 0.5510\% & 0.5295\% \\
		
		0.6 & 0.5932\% & 0.4655\% \\
		
		0.7 & 0.6445\% & 0.4228 \%\\
		
		0.8 & 0.6907\% & 0.3999\% \\
		
		0.9 & 0.7275\% & 0.3899\% \\
		
		1.0 & 0.7547\% & 0.3875\% \\
		\hline
	\end{tabular}
\end{minipage}	
\label{1.1error}
\end{table}

\begin{table}[hbt!]
	\caption{Relative errors at \da{different time steps} with $\alpha=1.5$ for Case 1 in Example 1. Left: $H=1/20$ and $l=6.$ Right: $H=1/40$ and $l=8.$}		
	\begin{minipage}[c]{0.5\textwidth}
		\centering
		\begin{tabular}{c c c}
			\hline
			$t$ & $e^{(1)}(t)$ & $e^{(2)}(t)$ \\
			\hline
			0.1 & 4.6637\% & 1.0388\% \\
			
			0.2 & 5.0996\% & 1.5235\% \\
			
			0.3 & 4.2770\% & 2.2505\% \\
			
			0.4 & 3.2231\% & 2.6408\% \\
			
			0.5 & 2.4167\% & 2.2978\% \\
			
			0.6 & 1.9937\% & 1.9367\% \\
			
			0.7 & 1.8844\% & 1.8525\% \\
			
			0.8 & 1.9575\% & 1.7751\% \\
			
			0.9 & 2.1136\% & 1.6038\% \\
			
			1.0 & 2.2911\% & 1.4076\% \\
			\hline
		\end{tabular}
	\end{minipage}
	\begin{minipage}[c]{0.5\textwidth}
		\centering
		\begin{tabular}{c c c}
			\hline
			$t$ & $e^{(1)}(t)$ & $e^{(2)}(t)$ \\
			\hline
			0.1 & 2.8759\% & 2.2597\% \\
			
			0.2 & 3.3907\% & 1.8471\% \\
			
			0.3 & 2.6049\% & 1.3953\% \\
			
			0.4 & 1.5417\% & 1.7346\% \\
			
			0.5 & 0.6943\%& 1.4025\% \\
			
			0.6 & 0.3017\% & 1.0359\% \\
			
			0.7 & 0.3072\% & 1.0668\% \\
			
			0.8 & 0.3176\% & 1.0312\% \\
			
			0.9 & 0.3493\% & 0.8011\% \\
			
			1.0 & 0.4593\% & 0.5027\% \\
			\hline
		\end{tabular}
	\end{minipage}
\label{1.5error}		
\end{table}

\begin{table}[hbt!]
	\caption{Relative errors at \da{different time steps} with $\alpha=1.9$ for Case 1 in Example 1. Left: $H=1/20$ and $l=6.$ Right: $H=1/40$ and $l=8.$}		
	\begin{minipage}[c]{0.5\textwidth}
		\centering
		\begin{tabular}{c c c}
			\hline
			$t$ & $e^{(1)}(t)$ & $e^{(2)}(t)$ \\
			\hline
			0.1 & 3.2835\% & 1.2876\% \\
			
			0.2 & 4.5091\% & 0.9312\% \\
			
			0.3 & 5.3425\% & 1.6117\% \\
			
			0.4 & 5.4649\% & 1.8804\% \\
			
			0.5 & 4.9539\% & 2.5350\% \\
			
			0.6 & 4.0613\% & 3.5565\% \\
			
			0.7 & 3.0423\% & 4.2390\% \\
			
			0.8 & 2.1076\% & 4.0676\% \\
			
			0.9 & 1.4487\% & 3.1557\% \\
			
			1.0 & 1.1994\% & 2.5295\% \\
			\hline
		\end{tabular}
	\end{minipage}
	\begin{minipage}[c]{0.5\textwidth}
		\centering
		\begin{tabular}{c c c}
			\hline
			$t$ & $e^{(1)}(t)$ & $e^{(2)}(t)$ \\
			\hline
			0.1 & 1.4581\% & 0.7435\% \\
			
			0.2 & 2.7035\% &  2.3210\% \\
			
			0.3 & 3.5888\% & 2.8675\% \\
			
			0.4 & 3.7700\% & 2.4045\% \\
			
			0.5 & 3.3087\% & 2.0839\% \\
			
			0.6 & 2.4581\% & 2.7680\% \\
			
			0.7 & 1.4846\% & 3.5367\% \\
			
			0.8 & 0.6802\% & 3.4962\% \\
			
			0.9 & 0.7163\% & 2.5995\% \\
			
			1.0 & 1.1724\% & 1.9349\% \\
			\hline
		\end{tabular}
	\end{minipage}
\label{1.9error}	
\end{table}

\subsubsection{Case 2: \da{M}ixed \da{time derivatives}}

\da{Let us now consider the case with} mixed time-fractional derivative \da{orders. In our numerical experiments, we set} $\alpha_1=1.1$ \da{and} $\alpha_2=1.9$. \da{Therefore, the solution field propagation is closer to the diffusion process in $\Omega_1$ (low-conductive subregion) and closer to the wave process in $\Omega_2$ (high-conductive subregion). In each continuum's equation of the resulting multicontinuum model, there is a linear combination of all the time derivatives. Figure \ref{m-1.1_1,9_f} presents distributions of the fine-grid solution field at different time steps. Again, we see the influence of the heterogeneous structure of the diffusion coefficient on the solution. However, the obtained distributions are more concentrated and less diffused due to the mixed time-fractional derivative orders.}

\begin{figure}[hbt!]
	\centering
	\includegraphics[width=0.29\textwidth]{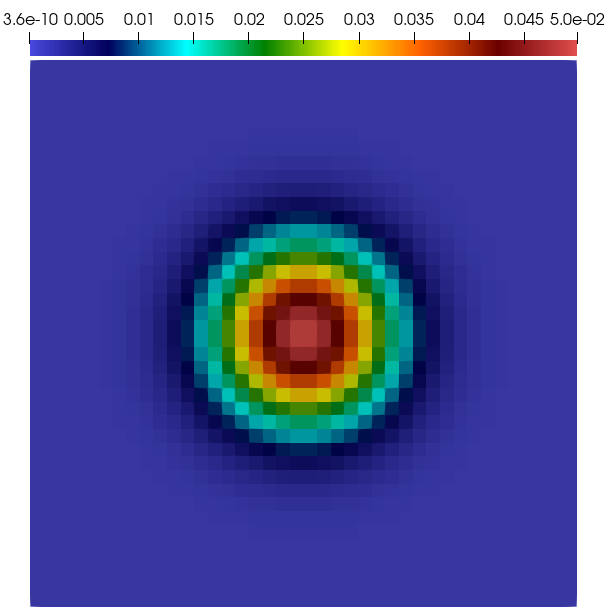}
	\includegraphics[width=0.29\textwidth]{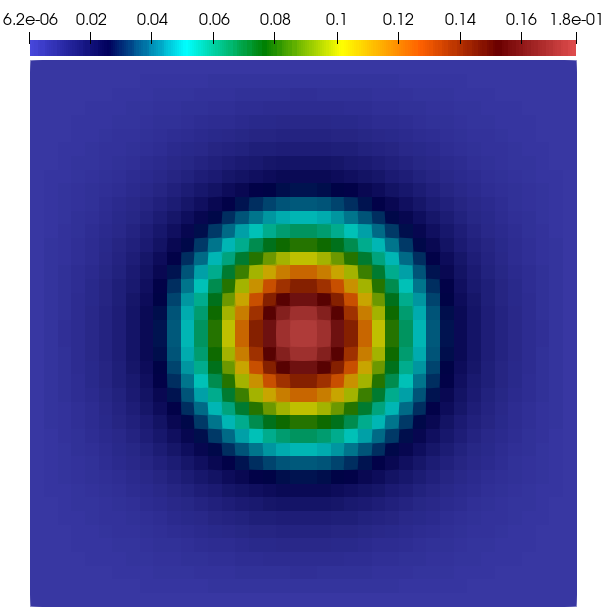}
	\includegraphics[width=0.29\textwidth]{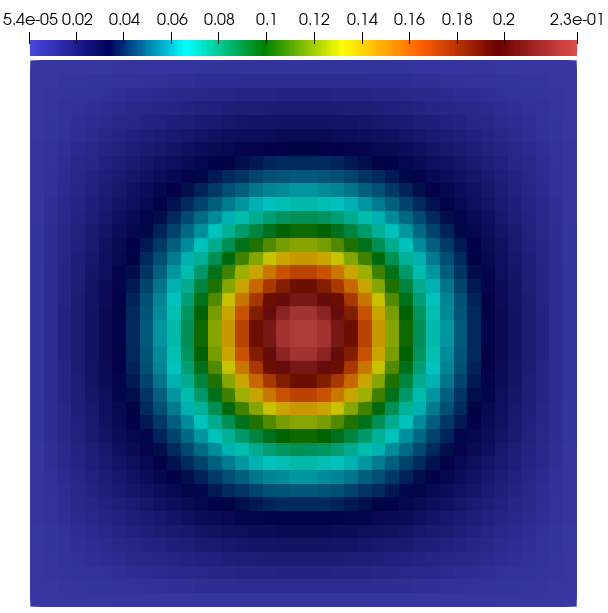}\\
	\includegraphics[width=0.29\textwidth]{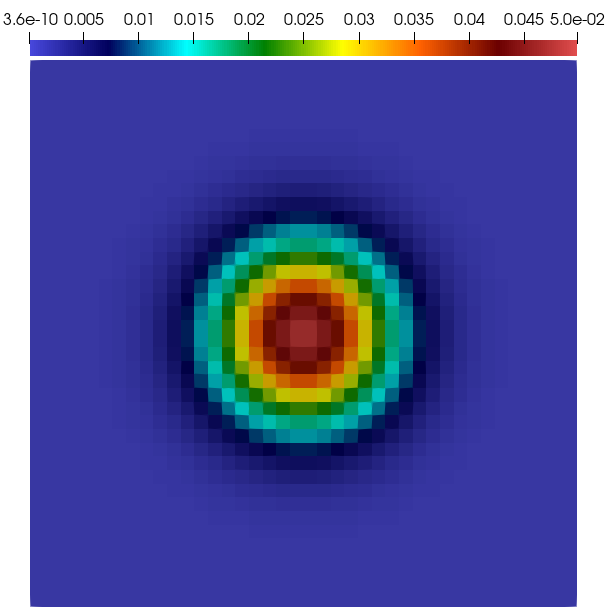}
	\includegraphics[width=0.29\textwidth]{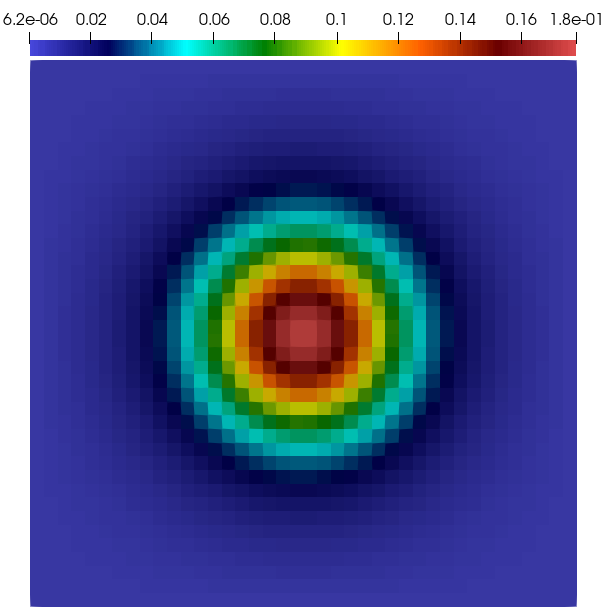}
	\includegraphics[width=0.29\textwidth]{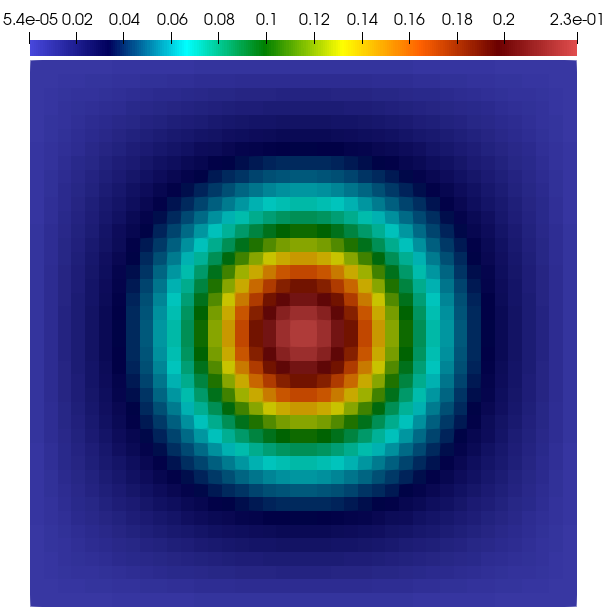}
	\caption{Average solution with $\alpha_1 = 1.1$, $\alpha_2 = 1.9$ at $t=0.1,0.5,1$ for Case 2 in Example 1 (from left to right)\da{, $H = 1/40$}. First row: reference averaged solution in $\Omega_1$. Second row: multiscale solution in $\Omega_1$.}
	\label{m-1.1-1.9-1}
\end{figure}

\begin{figure}[hbt!]
	\centering
	\includegraphics[width=0.29\textwidth]{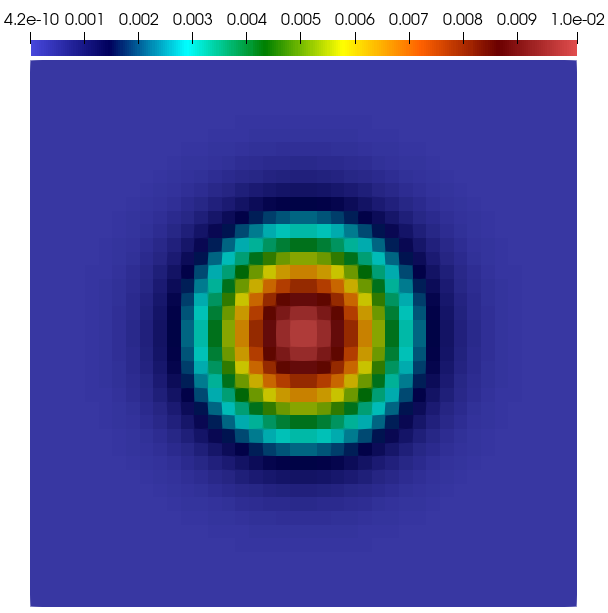}
	\includegraphics[width=0.29\textwidth]{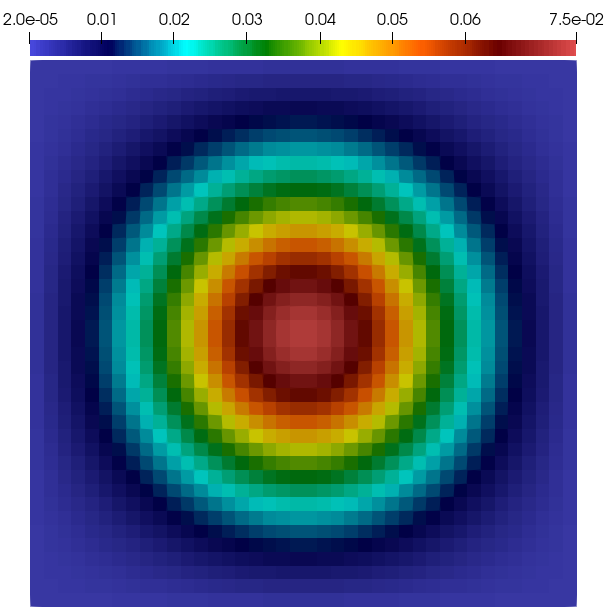}
	\includegraphics[width=0.29\textwidth]{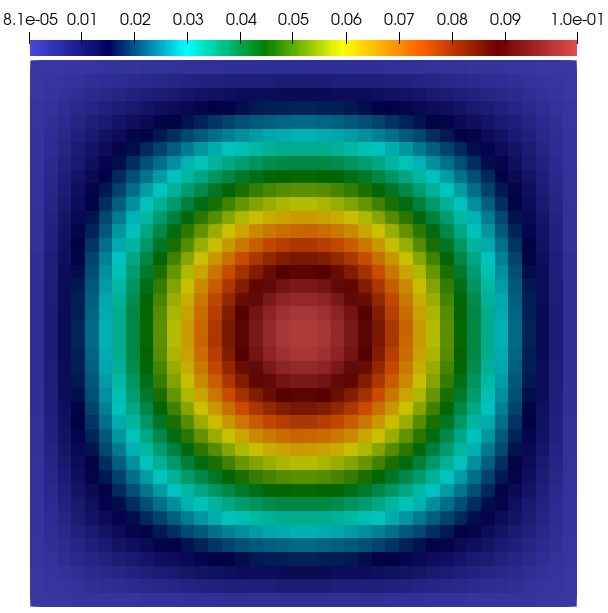}\\
	\includegraphics[width=0.29\textwidth]{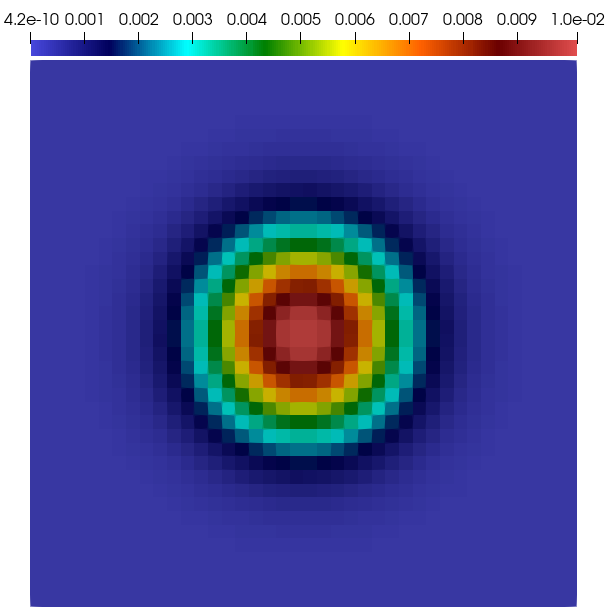}
	\includegraphics[width=0.29\textwidth]{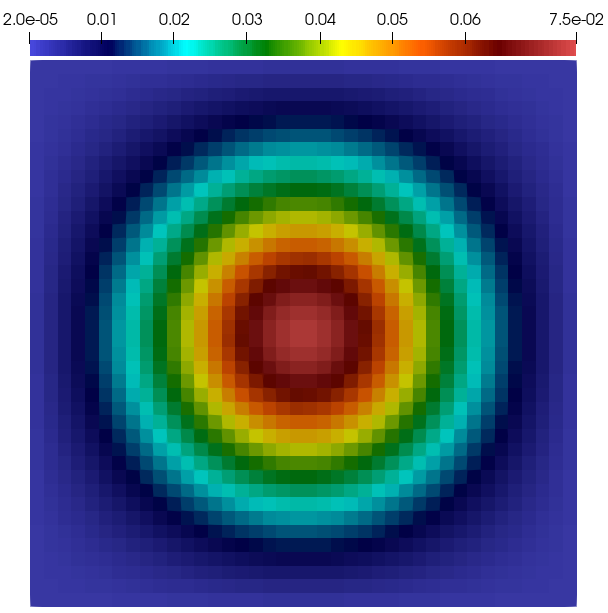}
	\includegraphics[width=0.29\textwidth]{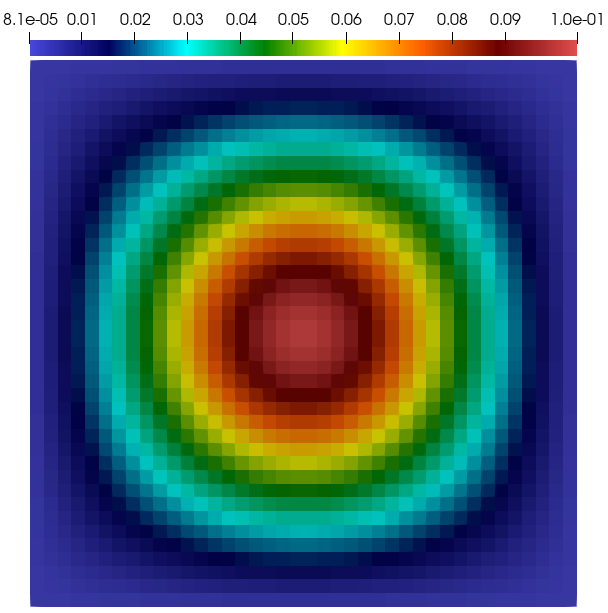}
	\caption{Average solution with $\alpha_1 = 1.1$, $\alpha_2 = 1.9$ at $t=0.1,0.5,1$ for Case 2 in Example 1 (from left to right)\da{, $H = 1/40$}. First row: reference averaged solution in $\Omega_2$. Second row: multiscale solution in $\Omega_2$.}
	\label{m-1.1-1.9-2}
\end{figure}

\da{Figures \ref{m-1.1-1.9-1} and \ref{m-1.1-1.9-2} present distributions of average solutions $H = 1/40$ in $\Omega_1$ and $\Omega_2$, respectively. In both figures, we depict the distributions of the reference and multiscale solutions from top to bottom. The obtained results are very similar, indicating the high accuracy of the proposed multicontinuum model.}

\begin{figure}[hbt!]
	\centering
	\includegraphics[width=0.29\textwidth]{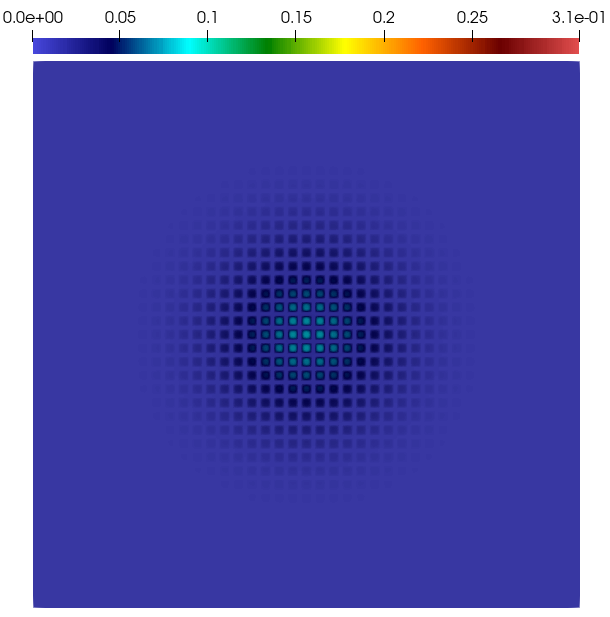}
	\includegraphics[width=0.29\textwidth]{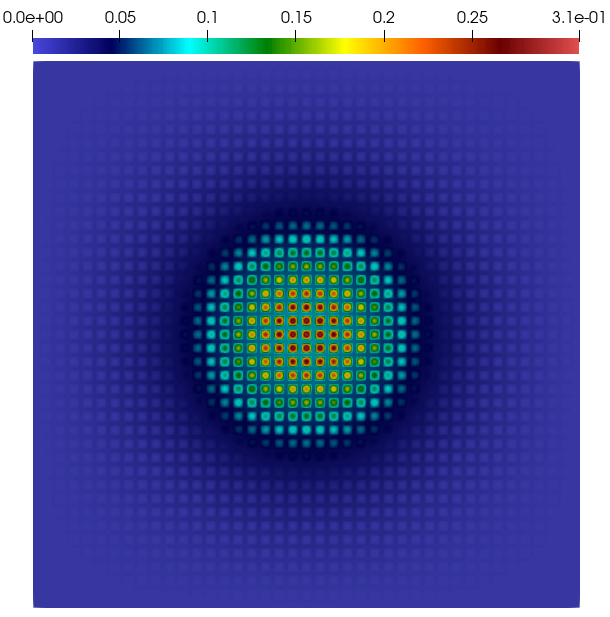}
	\includegraphics[width=0.29\textwidth]{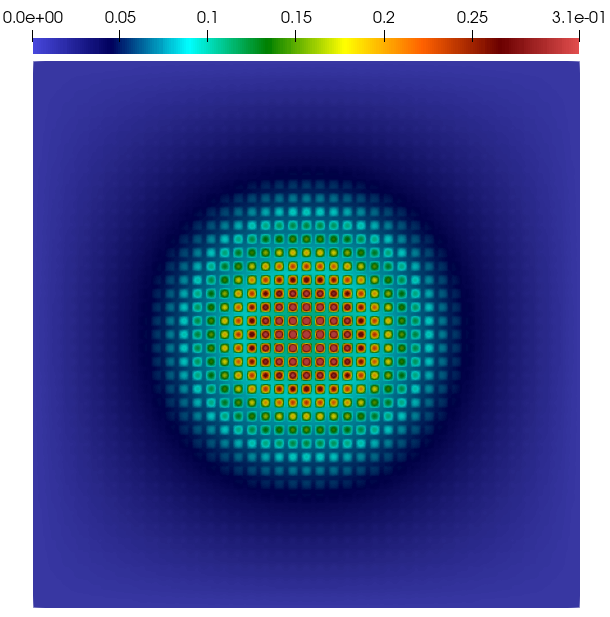}
	\caption{The fine-grid reference solution with $\alpha_1 = 1.1$, $\alpha_2 = 1.9$ at $t=0.1,0.5,1$ for Case 2 in Example 1 (from left to right).}
	\label{m-1.1_1,9_f}
\end{figure}


\da{Next, let us consider the errors of the multiscale solutions. Table \ref{c-1.1-1.9-error} presents the relative $L^2$ errors at different time steps for $\alpha_1 = 1.1$ and $\alpha_2 = 1.9$. From left to right, we depict the errors for the coarse grid sizes $H = 1/20$ and $H = 1/40$. One can see that the errors are small for both continua and coarse-grid sizes. We observe the errors decreasing with the decreasing of the coarse-grid size $H$, as expected.}

\begin{table}[hbt!]
	\caption{Relative errors at \da{different time steps} with $\alpha_1=1.1, \alpha_2=1.9$ for Case 2 in Example 1. Left: $H=1/20$ and $l=6.$ Right: $H=1/40$ and $l=8.$}		
	\begin{minipage}[c]{0.5\textwidth}
		\centering
		\begin{tabular}{c c c}
			\hline
			$t$ & $e^{(1)}(t)$ & $e^{(2)}(t)$ \\
			\hline
			0.1 & 5.0774\%  & 3.0555\%  \\
			
			0.2 & 3.0151\%  & 2.5509\%  \\
			
			0.3 & 2.8921\%  & 2.6385\%  \\
			
			0.4 & 2.8242\%  & 2.2639\%  \\
			
			0.5 & 2.6937\%  & 1.8062\%  \\
			
			0.6 & 2.5482\%  & 1.4270\% \\
			
			0.7 & 2.4243\%  & 1.1835\% \\
			
			0.8 & 2.3341\%  & 1.0884\% \\
			
			0.9 & 2.2852\%  & 1.1084\% \\
			
			1.0 & 2.2825\%  & 1.2035\% \\
			\hline
		\end{tabular}
	\end{minipage}
	\begin{minipage}[c]{0.5\textwidth}
		\centering
		\begin{tabular}{c c c}
			\hline
			$t$ & $e^{(1)}(t)$ & $e^{(2)}(t)$ \\
			\hline
			0.1 &  3.4044\%&  3.4382\%\\
			
			0.2 &  0.9572\% &  0.9120\%\\
			
			0.3 &  0.8672\%&  0.9831\%\\
			
			0.4 &  0.7755\%&  0.8696\%\\
			
			0.5 &  0.7030\%& 0.6527\%\\
			
			0.6 &  0.6655\%& 0.4656\%\\
			
			0.7 &  0.6574\%& 0.3468\%\\
			
			0.8 &  0.6573\%&  0.2877\%\\
			
			0.9 &  0.6604\%&  0.2733\%\\
			
			1.0 &  0.6756\%&  0.2944\%\\
			\hline
		\end{tabular}
	\end{minipage}
	\label{c-1.1-1.9-error}		
\end{table}

\FloatBarrier
\subsection{Example 2: Layered field}
\begin{figure}[hbt!]
	\centering
	\includegraphics[width=0.29\textwidth]{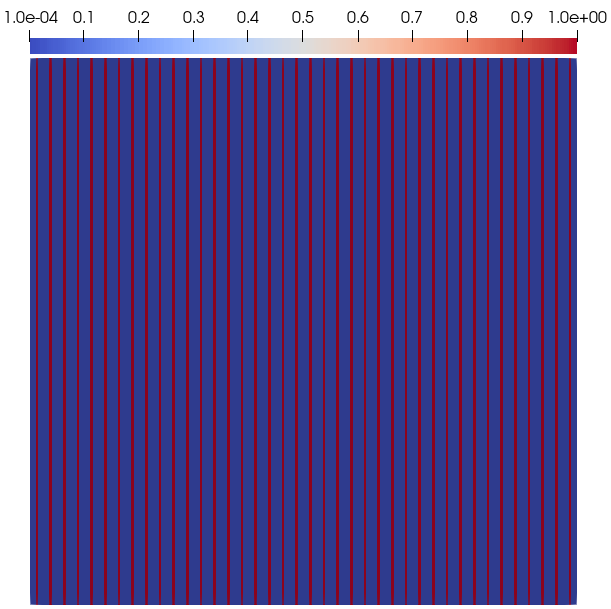}
	\caption{The coefficient $\kappa$ ($\Omega_1$: blue regions; $\Omega_2$: red regions). \da{Layered field.}}
	\label{kappa1}
\end{figure}

\da{As in the previous example, we define $\kappa$ using \eqref{eq:kappa}, but we consider a layered structure depicted in Figure \ref{kappa1}.} The source term $f$ is given as $f(x)=e^{-40((x_1 -0.5)^2+(x_2 -0.5)^2)}$ for any $x=(x_1,x_2)\in \Omega$. Also, we set the initial condition $u_0=u(x,0)=0.5e^{-40((x_1 -0.5)^2+(x_2 -0.5)^2)}$ and impose zero Neumann boundary condition 
\begin{equation}
	\frac{\partial u}{\partial n}=0, \quad x \in \partial \Omega.
\end{equation}

\da{Again, we consider both regular and mixed cases of time-fractional derivative orders in the following numerical experiments.}

\subsubsection{Case 1: \da{R}egular \da{time derivatives}}

\begin{figure}[hbt!]
	\centering
	\includegraphics[width=0.29\textwidth]{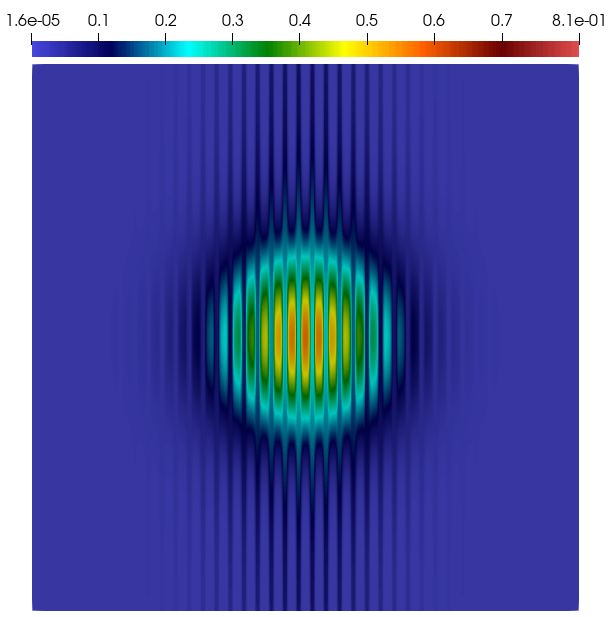}
	\includegraphics[width=0.29\textwidth]{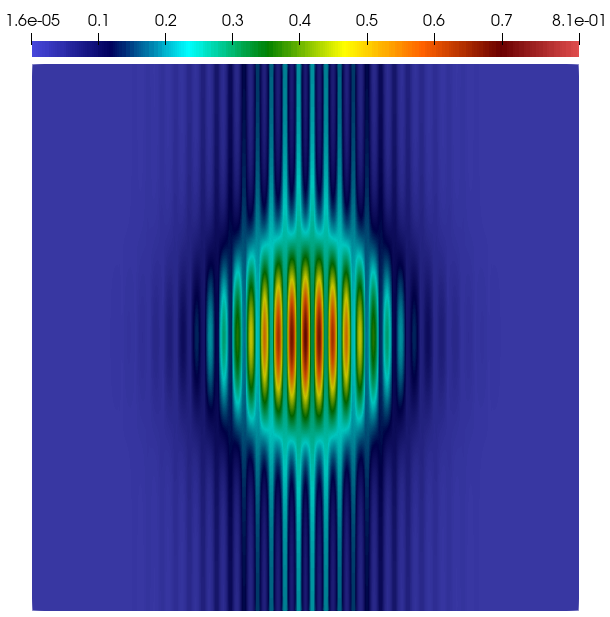}
	\includegraphics[width=0.29\textwidth]{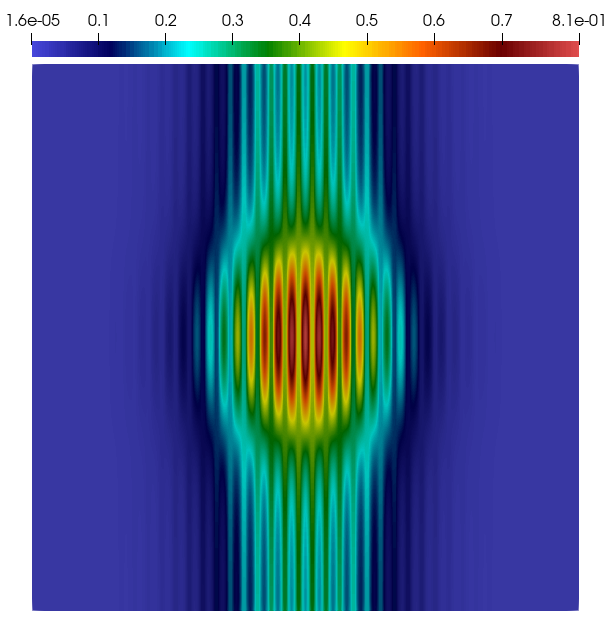}
	\caption{The fine-grid reference solution with $\alpha = 1.2$ at $t=0.1,0.5,1$ for Case 1 in Example 2 (from left to right).}
	\label{1.2_f}
\end{figure}

In this case, we set the time fractional order derivative $\alpha=1.2$ \da{for both continua}. \da{Figure \ref{1.2_f} depicts distributions of the fine-grid solution field at different time steps. One can see the influence of the layered structure of the heterogeneous diffusion coefficient. The initial distribution of the solution gradually propagates along the vertical directions due to the high conductivity of the channels. The obtained numerical results correspond to the simulated process.}

\begin{figure}[hbt!]
	\centering
	\includegraphics[width=0.29\textwidth]{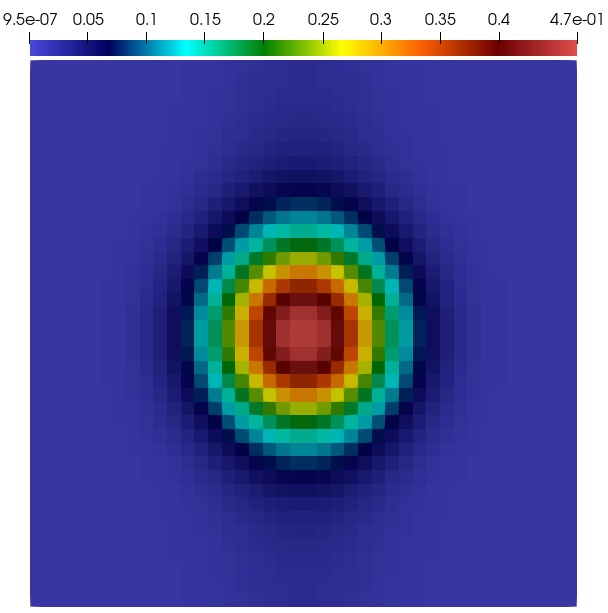}
	\includegraphics[width=0.29\textwidth]{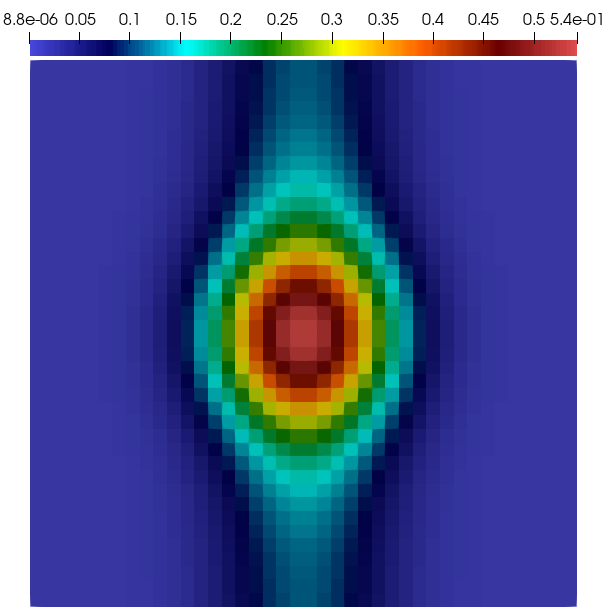}
	\includegraphics[width=0.29\textwidth]{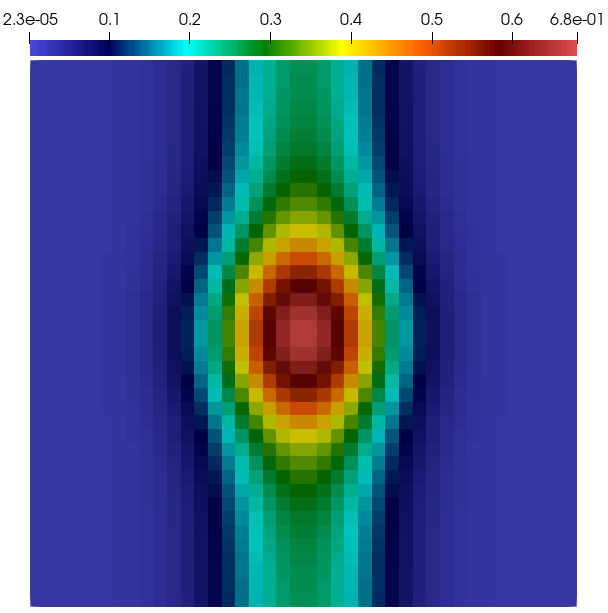}\\
	\includegraphics[width=0.29\textwidth]{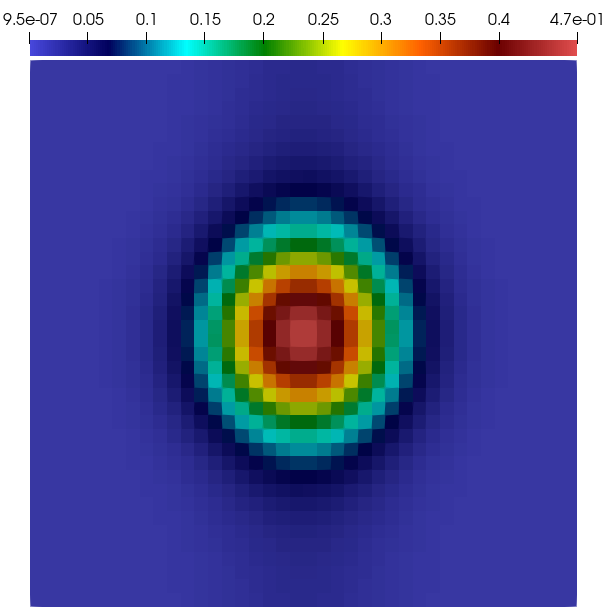}
	\includegraphics[width=0.29\textwidth]{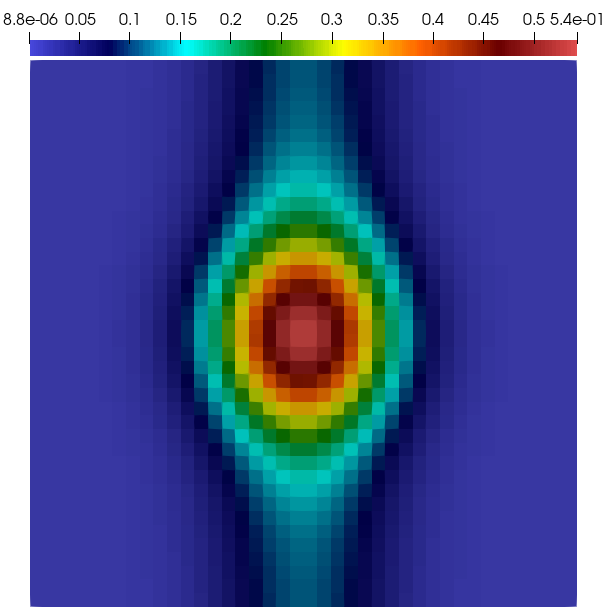}
	\includegraphics[width=0.29\textwidth]{l-1.2-f-1-1.png}
	\caption{Average solution with $\alpha = 1.2$ at $t=0.1,0.5,1$ for Case 1 in Example 2 (from left to right)\da{, $H = 1/40$}. First row: reference averaged solution in $\Omega_1$. Second row: multiscale solution in $\Omega_1$.}
	\label{1.2_1}
\end{figure}

\begin{figure}[hbt!]
	\centering
	\includegraphics[width=0.29\textwidth]{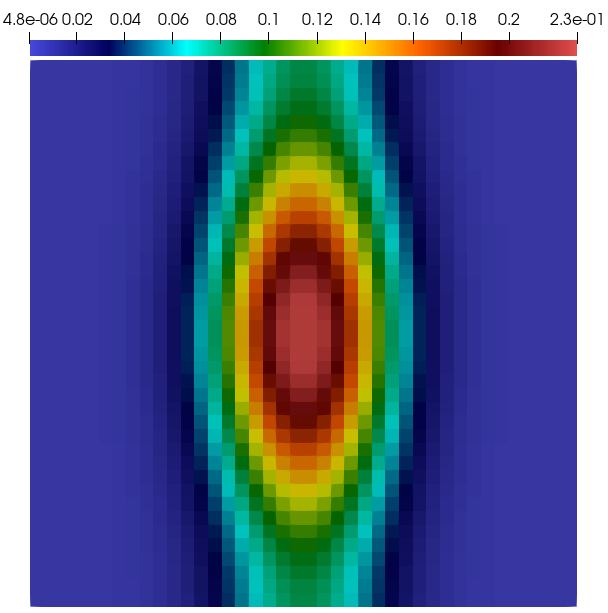}
	\includegraphics[width=0.29\textwidth]{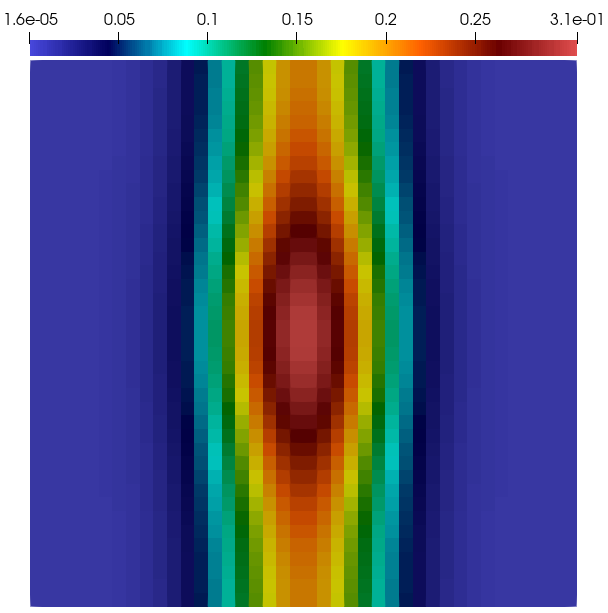}
	\includegraphics[width=0.29\textwidth]{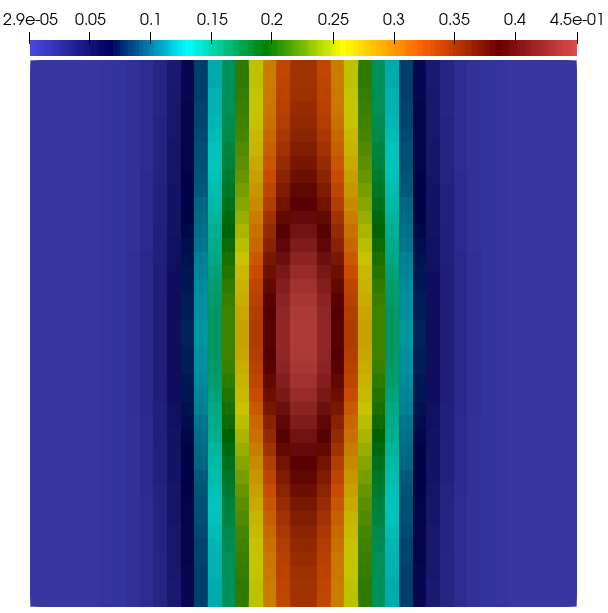}\\
	\includegraphics[width=0.29\textwidth]{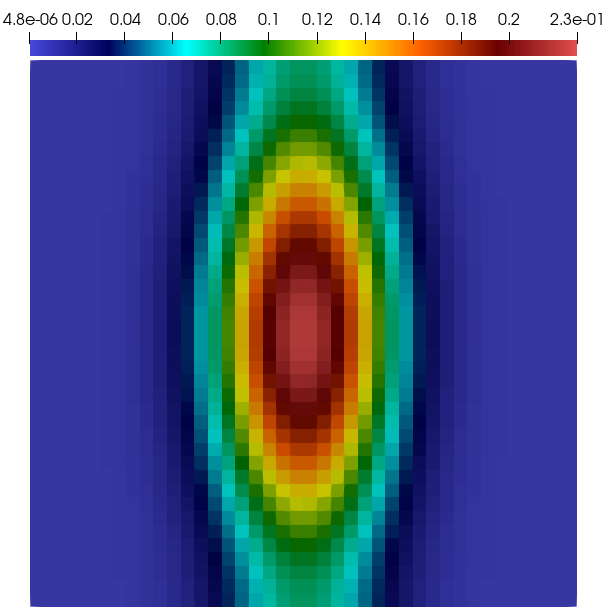}
	\includegraphics[width=0.29\textwidth]{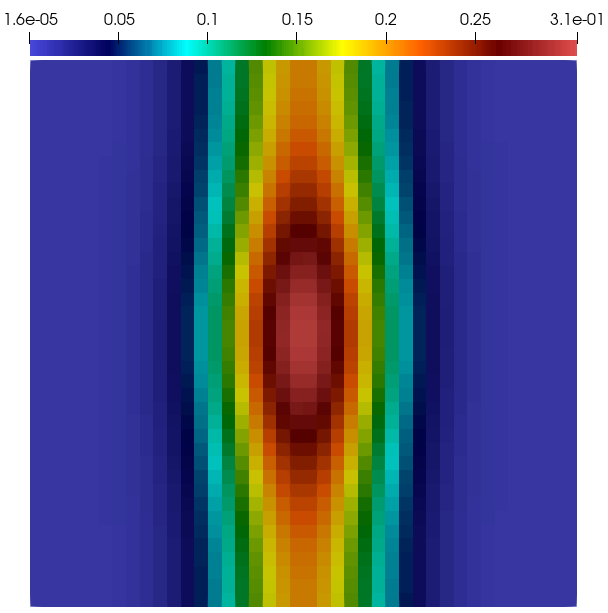}
	\includegraphics[width=0.29\textwidth]{l-1.2-f-2-1.png}
	\caption{Average solution with $\alpha = 1.2$ at $t=0.1,0.5,1$ for Case 1 in Example 2 (from left to right)\da{, $H = 1/40$}. First row: reference averaged solution in $\Omega_2$. Second row: multiscale solution $U_2$ in $\Omega_2$.}
	\label{1.2_2}
\end{figure}

\da{Next, let us consider the average solutions of the reference and multiscale solutions. In Figures \ref{1.2_1} and \ref{1.2_2}, we present distributions of the average solutions in $\Omega_1$ and $\Omega_2$ with $H = 1/40$, respectively. From top to bottom, we depict the reference and multiscale solutions. As in the previous example, we see that the obtained results are very similar. Moreover, one can notice that the solution propagation has an anisotropic nature due to the structure of the high-contrast diffusion coefficient. We observe more rapid propagation of the average solution in $\Omega_2$ than in $\Omega_1$, as expected.}


\da{Let us consider the errors of the multiscale solution. In Table \ref{1.2error}, we present the relative $L^2$ errors of the average multiscale solutions at different time steps. From left to right, we have the errors for coarse-grid sizes $H=1/20$ and $H=1/40$. One can see that the errors are minor for both coarse grids. One can see that the errors gradually decrease with time. Moreover, the smaller the coarse-grid size, the higher the accuracy we achieve. Therefore, our multicontinuum approach can approximate the reference solution with high accuracy.}

\begin{table}[hbt!]
	\caption{Relative errors at \da{different time steps} with $\alpha=1.2$ for Case 1 in Example 2. Left: $H=1/20$ and $l=6.$ Right: $H=1/40$ and $l=8.$}		
	\begin{minipage}[c]{0.5\textwidth}
		\centering
		\begin{tabular}{c c c}
			\hline
			$t$ & $e^{(1)}(t)$ & $e^{(2)}(t)$ \\
			\hline
			0.1 & 3.4643\%  & 3.0195\%  \\
			
			0.2 & 2.3059\%  & 3.7543\%  \\
			
			0.3 & 2.0710\%  & 1.9225\%  \\
			
			0.4 & 2.1109\%  & 1.5818\%  \\
			
			0.5 & 2.1660\%  & 1.5747\%  \\
			
			0.6 & 2.1843\%  & 1.5947\% \\
			
			0.7 & 2.1728\%  & 1.5984\% \\
			
			0.8 & 2.1416\%  & 1.5965\% \\
			
			0.9 & 2.1021\%  & 1.5849\% \\
			
			1.0 & 2.0577\%  & 1.5777\% \\
			\hline
		\end{tabular}
	\end{minipage}
	\begin{minipage}[c]{0.5\textwidth}
		\centering
		\begin{tabular}{c c c}
			\hline
			$t$ & $e^{(1)}(t)$ & $e^{(2)}(t)$ \\
			\hline
			0.1 & 1.7997\%   & 1.9799\% \\
			
			0.2 & 0.5696\%  &  2.8582\%\\
			
			0.3 &  0.4282\%& 0.8036\% \\
			
			0.4 &  0.4650\%&  0.4468\%\\
			
			0.5 & 0.5128\% & 0.4341\%\\
			
			0.6 & 0.5614\% & 0.4056\%\\
			
			0.7 & 0.5940\% & 0.4081\%\\
			
			0.8 & 0.6084\% &  0.4162\%\\
			
			0.9 & 0.6124\% &  0.4143\%\\
			
			1.0 & 0.6078\% &  0.4154\%\\
			\hline
		\end{tabular}
	\end{minipage}
	\label{1.2error}		
\end{table}

\subsubsection{Case 2: \da{M}ixed \da{time derivatives}}

\begin{figure}[hbt!]
	\centering
	\includegraphics[width=0.29\textwidth]{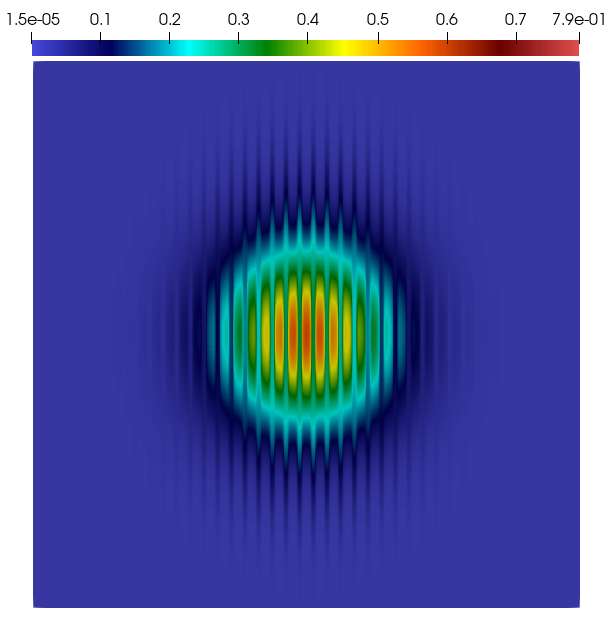}
	\includegraphics[width=0.29\textwidth]{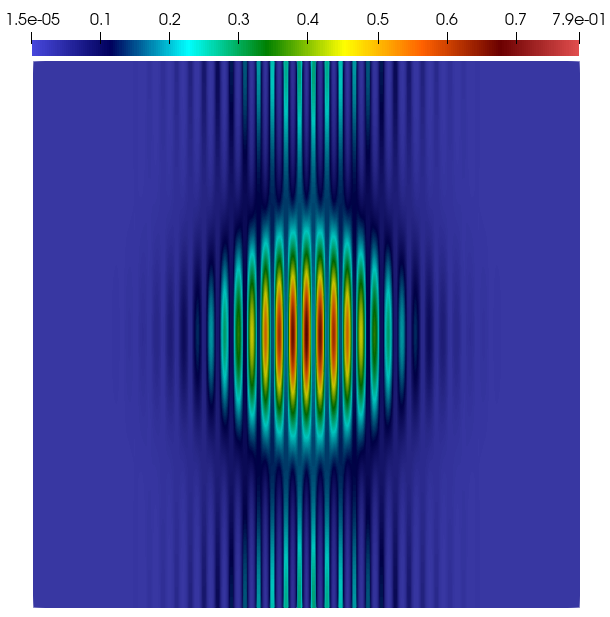}
	\includegraphics[width=0.29\textwidth]{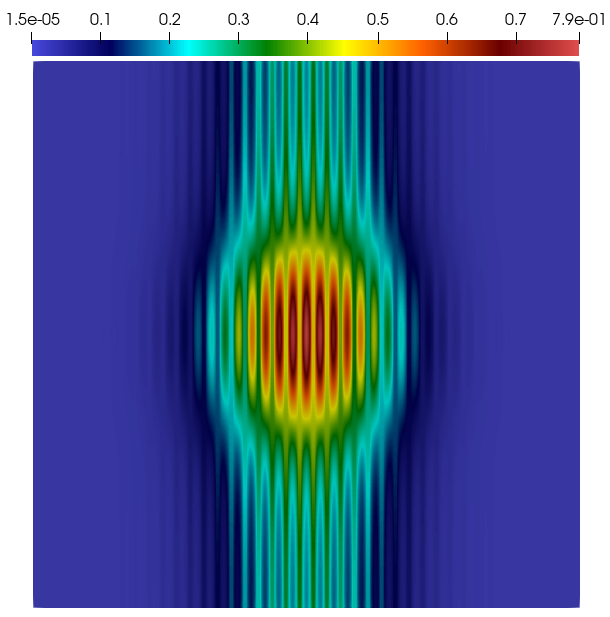}
	\caption{The fine-grid reference solution with $\alpha_1 = 1.1$, $\alpha_2 = 1.9$ at $t=0.1,0.4,1$ for Case 2 in Example 2 (from left to right).}
	\label{l-m-f}
\end{figure}

\da{Finally, let us consider the case with mixed time-fractional derivative orders. Again, we set $\alpha_1 = 1.1$ in $\Omega_1$ and $\alpha_2 = 1.9$ in $\Omega_2$. Therefore, the solution propagation in the high-conductive channels is closer to the wave process, while the spread in the low-conductive subregions is closer to the diffusion process. Figure \ref{l-m-f} shows distributions of the fine-grid solution at different time steps. One can see the changes in the solution distributions compared to the previous case due to the mixed time-fractional derivatives. It is especially noticeable in the high-conductive channels at the time instant $t = 0.4$.}

\begin{figure}[hbt!]
	\centering
	\includegraphics[width=0.29\textwidth]{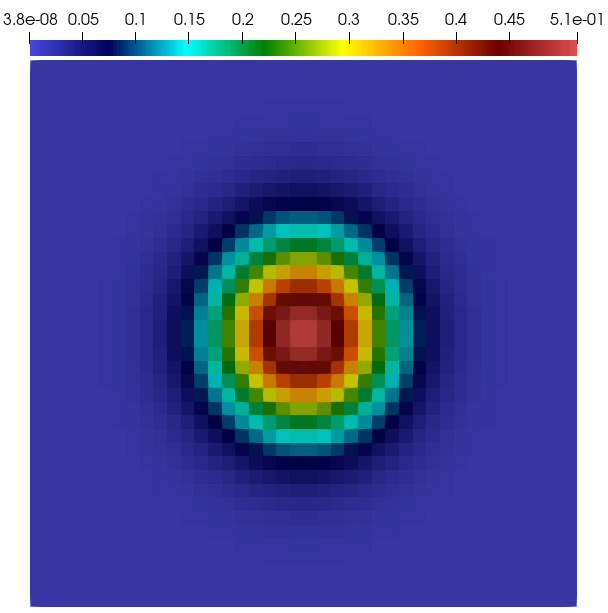}
	\includegraphics[width=0.29\textwidth]{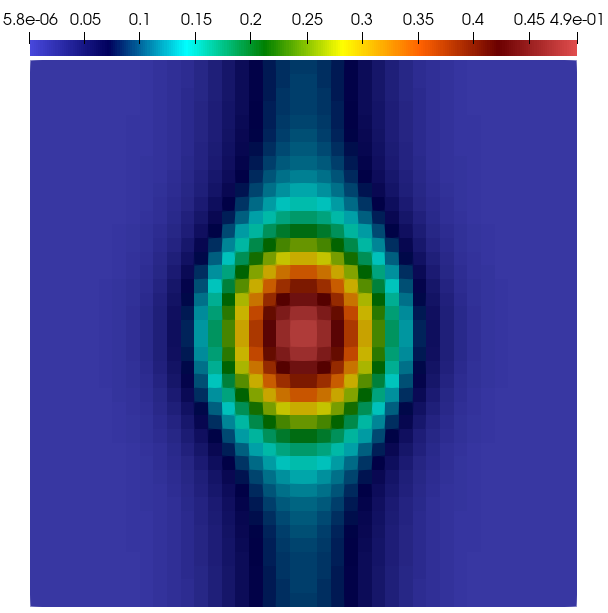}
	\includegraphics[width=0.29\textwidth]{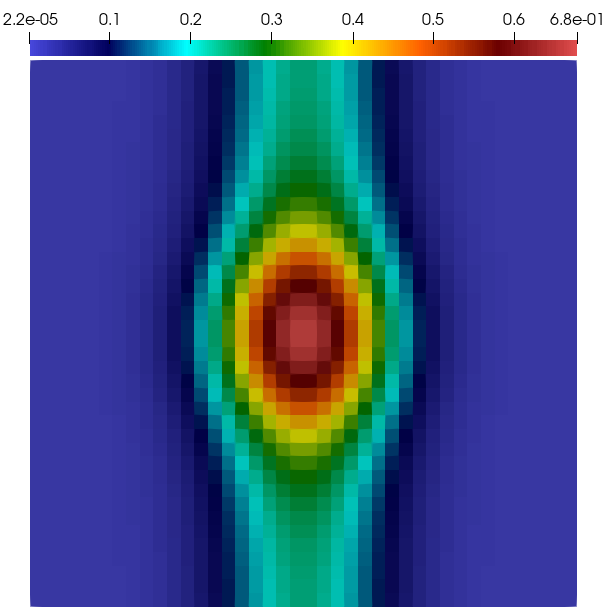}\\
	\includegraphics[width=0.29\textwidth]{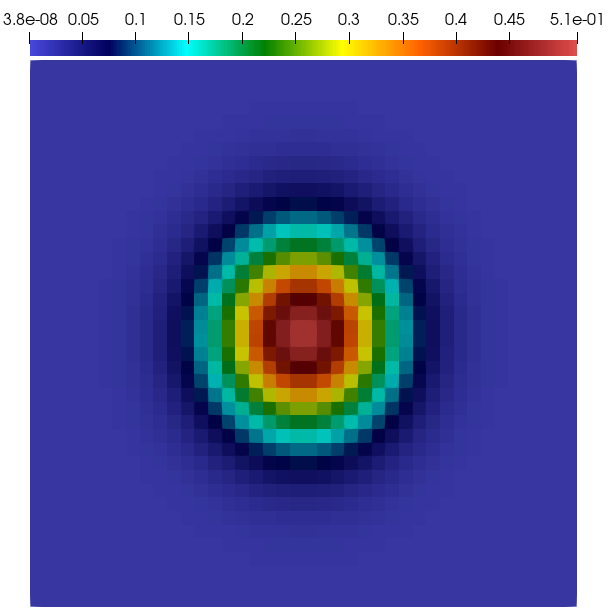}
	\includegraphics[width=0.29\textwidth]{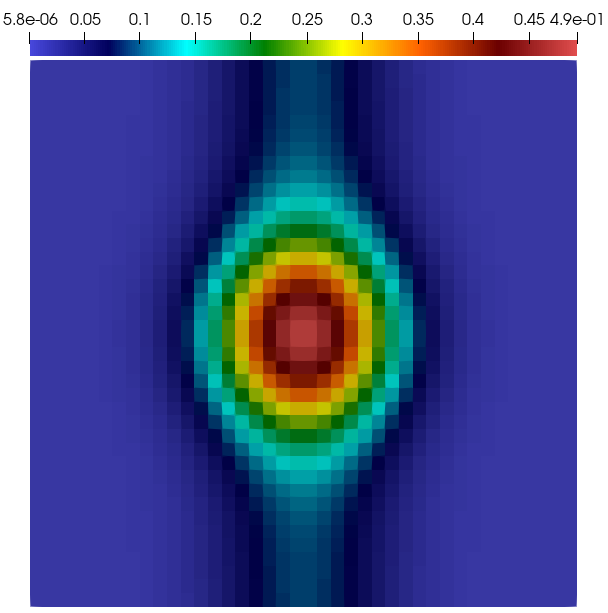}
	\includegraphics[width=0.29\textwidth]{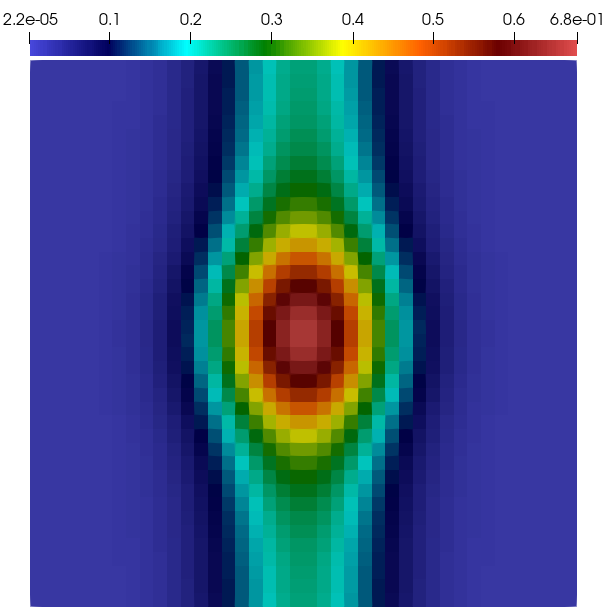}
	\caption{Average solution with $\alpha_1 = 1.1$, $\alpha_2 = 1.9$ at $t=0.1,0.4,1$ for Case 2 in Example 2 (from left to right)\da{, $H = 1/40$}. First row: reference averaged solution in $\Omega_1$. Second row: multiscale solution in $\Omega_1$.}
	\label{1.1-1.9-1}
\end{figure}

\begin{figure}[hbt!]
	\centering
	\includegraphics[width=0.29\textwidth]{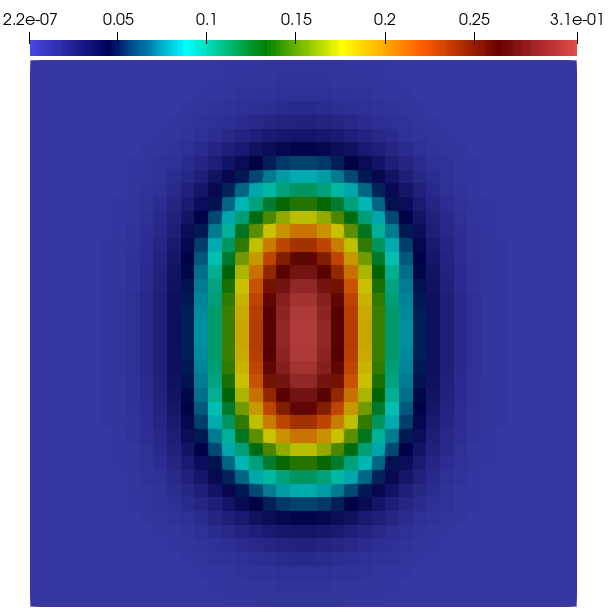}
	\includegraphics[width=0.29\textwidth]{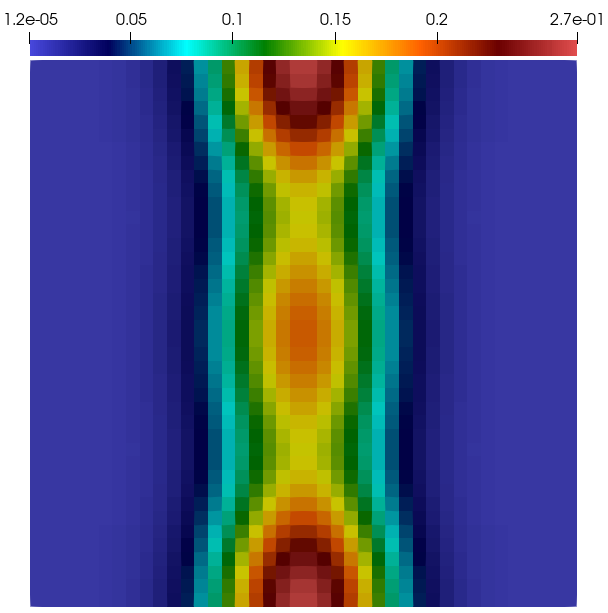}
	\includegraphics[width=0.29\textwidth]{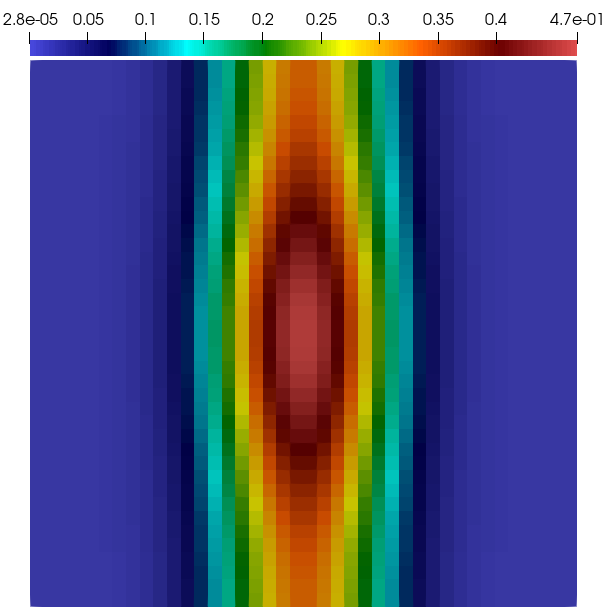}\\
	\includegraphics[width=0.29\textwidth]{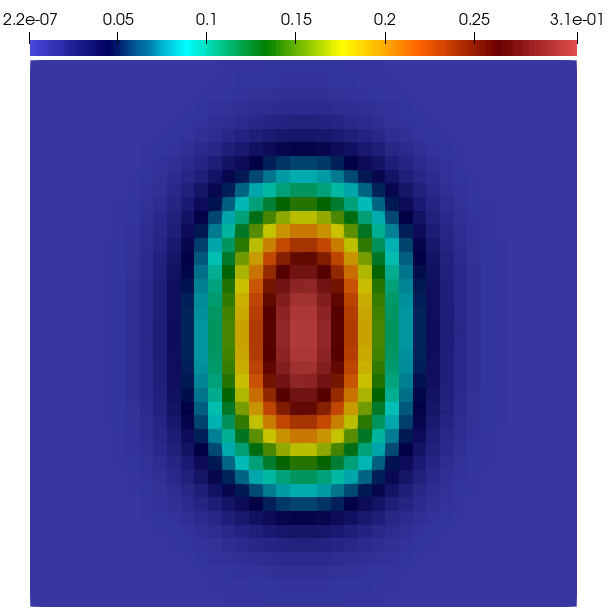}
	\includegraphics[width=0.29\textwidth]{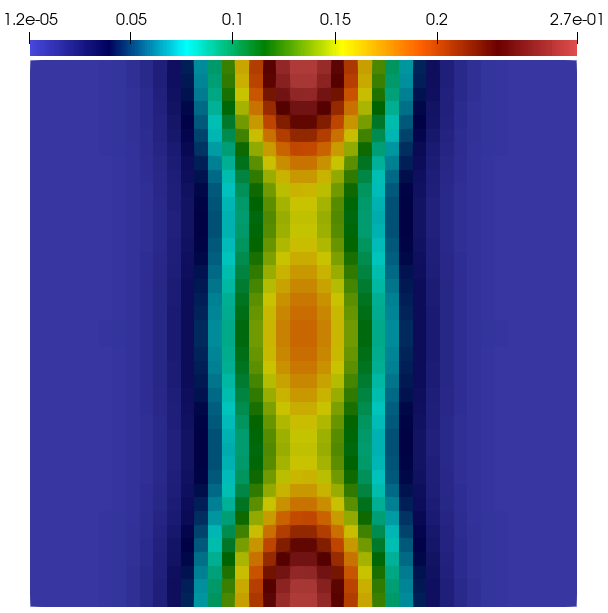}
	\includegraphics[width=0.29\textwidth]{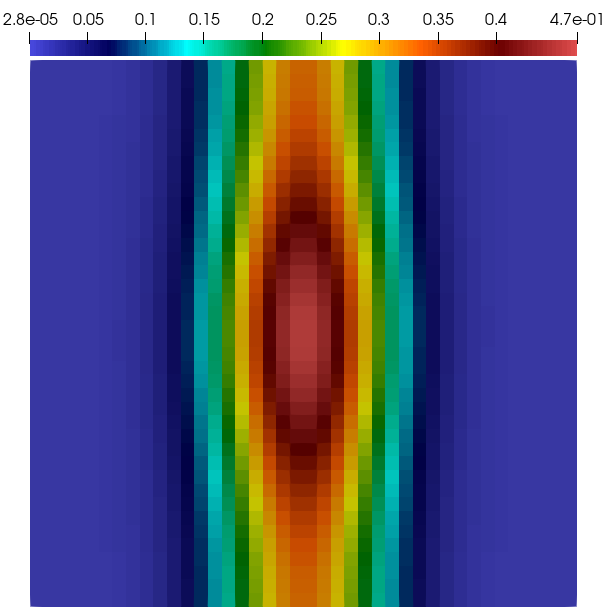}
	\caption{Average solution with $\alpha_1 = 1.1$, $\alpha_2 = 1.9$ at $t=0.1,0.4,1$ for Case 2 in Example 2 (from left to right)\da{, $H = 1/40$}. First row: reference averaged solution in $\Omega_2$. Second row: multiscale solution in $\Omega_2$.}
	\label{1.1-1.9-2}
\end{figure}

\da{Figures \ref{1.1-1.9-1} and \ref{1.1-1.9-2} present distributions of the average solution fields with $H=1/40$ in $\Omega_1$ and $\Omega_2$, respectively. We depict the reference average solution on the top and the multiscale solution on the bottom. One can see that our multiscale solutions are very similar to the reference ones, indicating high accuracy. Compared to the previous case, the average solution in $\Omega_1$ changed regarding ranges but mostly remained the dynamics. In contrast, we can see the changes in the dynamics of the average solution in $\Omega_2$ due to the mixed time-fractional derivatives. The differences are more noticeable at the time instants $t = 0.1$ and $t = 0.4$.}


\da{Next, let us consider the errors of our multicontinuum approach. Table \ref{1.1-1.9-error} presents the relative $L^2$ errors at different time steps for both continua using $H=1/20$ and $H=1/40$. One can see that all the errors are minor. We can observe the reduction in the errors with a decrease in the coarse-grid size $H$.}

\begin{table}[hbt!]
	\caption{Relative errors at \da{different time steps} with $\alpha_1=1.1, \alpha_2=1.9$ for Case 2 in Example 2. Left: $H=1/20$ and $l=6.$ Right: $H=1/40$ and $l=8.$}		
	\begin{minipage}[c]{0.5\textwidth}
		\centering
		\begin{tabular}{c c c}
			\hline
			$t$ & $e^{(1)}(t)$ & $e^{(2)}(t)$ \\
			\hline
			0.1 & 3.6809\%  & 2.6249\%  \\
			
			0.2 & 2.9735\%  & 2.8666\%  \\
			
			0.3 & 2.3455\%  & 2.8279\%  \\
			
			0.4 & 2.1554\%  & 2.1954\%  \\
			
			0.5 & 2.3374\%  & 2.2074\%  \\
			
			0.6 & 2.2823\%  & 1.8027\% \\
			
			0.7 & 2.0348\%  & 1.6265\% \\
			
			0.8 & 2.0175\%  & 1.5654\% \\
			
			0.9 & 2.1512\%  & 1.7002\% \\
			
			1.0 & 2.1917\%  & 1.7845\% \\
			\hline
		\end{tabular}
	\end{minipage}
	\begin{minipage}[c]{0.5\textwidth}
		\centering
		\begin{tabular}{c c c}
			\hline
			$t$ & $e^{(1)}(t)$ & $e^{(2)}(t)$ \\
			\hline
			0.1 &   1.8825\% & 0.6262\% \\
			
			0.2 &   1.1861\%&  0.7398\%\\
			
			0.3 &  0.7846\%&  0.7883\%\\
			
			0.4 &  0.5348\%&  0.9551\%\\
			
			0.5 &  0.7289\%& 1.2642\%\\
			
			0.6 &  0.6791\%& 0.9199\%\\
			
			0.7 &  0.5695\%& 0.5843\%\\
			
			0.8 &  0.6366\%&  0.6088\%\\
			
			0.9 &  0.7034\%&  0.7133\%\\
			
			1.0 &  0.6888\%&  0.5972\%\\
			\hline
		\end{tabular}
	\end{minipage}
	\label{1.1-1.9-error}		
\end{table}

\da{Therefore, our multicontinuum approach can provide high accuracy for different types of heterogeneities that can result in isotropic and anisotropic propagation. Moreover, it can accurately capture different time-fractional derivative orders, including the case with different orders in different continua.}

	\FloatBarrier
	\section{Conclusion}
	
	In this work, we \da{have} appl\da{ied} the multicontinuum homogenization method to \da{derive} the \da{multicontinuum} time-fractional diffusion-wave \da{model}. \da{First, we formulated constraint} cell problems \da{in oversampled regions} considering gradient effects and averages. By \da{solving the} constraint cell problems, we \da{obtained multicontinuum expansions of the fine-scale solution}. \da{Then, we rigorously derived the multicontinuum models for two cases of time-fractional derivatives. In the first case, time-fractional orders are the same for all continua. The second case represents the mixed time-fractional derivatives, where we have different orders in different continua. To check the effectiveness of the proposed multicontinuum approach, we have conducted numerical experiments by solving model problems with two continua. The results show that our multicontinuum models can provide high accuracy for different heterogeneous media and time-fractional derivatives.}

\section*{Acknowledgments}
	
The research of Huiran Bai is supported by the Postgraduate Scientific Research Innovation Project of Xiangtan University (XDCX2024Y160) and the Chinese Government Scholarship (CSC No. 202408430165). 

The research of Yin Yang and Wei Xie is supported by the National Natural Science Foundation of China Project (No. 12261131501), National Foreign Experts Program (No.S20240066), the Project of Scientific Research Fund of the Hunan Provincial Science and Technology Department (No.2023GK2029, No.2024JC1003, No.2024JJ1008), and “Algorithmic Research on Mathematical Common Fundamentals” Program for Science and Technology Innovative Research Team in Higher Educational Institutions of Hunan Province of China. 

The research of Dmitry Ammosov is supported by the Khalifa University Postdoctoral Fellowship Program.
	
	
	
	\bibliography{ref}	
\end{document}